\documentclass[12pt,a4paper,leqno]{article}
\usepackage[latin1]{inputenc}
\usepackage{amsmath}
\usepackage{amsfonts}
\usepackage{amssymb}
\usepackage{mathrsfs}
\usepackage{makeidx}
\usepackage{graphicx}
\usepackage{color}
\usepackage{multirow}
\usepackage[left=3.00cm, right=3.00cm, top=3.00cm, bottom=3.00cm]{geometry}
\usepackage{hyperref}
\hypersetup{ colorlinks=true
pdfborder=0 0 0}
\newcounter{subs}[section]
\renewcommand{\thesubs}{\textbf{\thesection.\arabic{subs}}}
\newcommand{\saa}{\refstepcounter{subs}\thesubs~}
\newcounter{claim}[section]

\numberwithin{equation}{section}

\newcommand{\ra}{\rightarrow}

\newcommand{\mb}[1]{\mathbf{#1}}
\newcommand{\mc}[1]{\mathcal{#1}}
\renewcommand{\subset}{\subseteq}
\newcommand {\qed} {\hfill $ \square $\vspace{.5cm}}

\newcommand{\caa}{\mb{c}}
\newcommand{\Hm}{\mathcal{H}}

\newcommand{\ed}{{\dot{=}}}
\newcommand{\edq}{{\dot{\equiv}}}
\newcommand{\x}[1]{\xi_{#1}}
\newcommand{\T}[1]{T_{#1}}
\allowdisplaybreaks
\newtheorem{thm}[subs]{Theorem}
\newtheorem{prop}[subs]{Proposition}
\newtheorem{lemma}[subs]{Lemma}
\newtheorem{cor}[subs]{Corollary}
\newtheorem{defn}[subs]{Definition}
\newtheorem{assumption}[subs]{Assumption}

\newtheorem{remark}[subs]{Remark}
\newtheorem{notation}[subs]{Notation}
\author{Xun Xie\footnote{Email:xieg7@163.com. Date:09/13/2015.}}
\title{A decomposition formula for the Kazhdan-Lusztig basis of affine Hecke algebras of rank 2}
\date{}

\begin{document}
\maketitle
\abstract{
In this paper, we prove a decomposition formula for the Kazhdan-Lusztig basis of affine Hecke algebras of rank 2 with positive weight function. Then we discuss some applications of this kind of  decomposition  to Lusztig's conjectures P1-P15.
}
\section{Introduction}
\saa For each Coxeter group, we have left cells, right cells, and two-sided cells, defined by  Kazhdan-Lusztig basis of the corresponding Hecke algebra, see\cite{kazhdan_lusztig79representation}.
The cells of an affine Weyl group are studied in \cite{lusztig1985cellsI,lusztig1987cellsII,lusztig1987III,lusztig1989cellsIV}. It is shown that the cells of affine Weyl groups have many remarkable properties and play an important role in representations of affine Hecke algebras. This  theory about two-sided cells  is  also discussed by \cite{lusztig2003hecke} in the context of  Coxeter groups with unequal parameters. Some important properties of two-sided cells are contained in the conjectures P1-P15, see \cite[\S14]{lusztig2003hecke}. In general, these conjectures  are still open   for   Coxeter groups with unequal parameters, and even for finite Weyl group. One of our motivation is to try to understand these conjectures for affine Weyl groups of rank 2. However, we do not obtain a complete proof of P1-P15 in this case. The main result  of this paper is a kind of decomposition formula for Kazhdan-Lusztig basis (Theorem \ref{th:dec}), which will reduce P1-P15 in our case to proving P8 and determining Lusztig's $ \mb{a} $-function explicitly. 

Conjectures P1-P15 have been proved in the case of  Coxeter groups with equal parameter, as a corollary of the positivity conjecture, see \cite{ben_will2014hodge} and \cite{lusztig2014unequal}. In this paper, unless otherwise specified, P1-P15 refer to these conjectures for unequal parameters. It is known that 
positivity does not hold  for Coxeter groups with unequal parameters, which makes P1-P15 still open. 
For finite Coxeter groups that admit unequal parameters (type $B_n(n\geq3)$, $I(m)$ ($m$ even) and $F_4$), we have already known that P1-P15 hold for type $I(m)$, $F_4$ and for type $B_n$ with ``asymptotic parameters", see  \cite[Thm.5.3]{geck2011iso} and the references therein.

For an affine Weyl group, the author has proved that   P1-P15 hold on  a distinguished two-sided cell---the lowest two-sided cell, see \cite{xie2015lowest}.
The main tool for this   is a  decomposition formula of Kazhdan-Lusztig basis related to the lowest two-sided cell. To generalize the method used in \cite{xie2015lowest}, we are motivated to investigate the case of  affine Weyl groups of rank 2  and  focus on type $ \tilde{B}_2 $ and $ \tilde{G} _2$. The partitions of two-sided cells of an affine Weyl group of type $ \tilde{B}_2 $ and $ \tilde{G} _2$ for  various parameters have been  worked out in \cite{guilhot0810.5165,guilhot2010rank2}, which provide us with a lot of nontrivial examples of cells and make it possible to check some properties directly by calculation. This is the reason why we are concerned with  affine Weyl groups of rank 2. As the case of the lowest two-sided cell, our first step is to establish a decomposition formula for Kazhdan-Lusztig basis of affine Weyl groups of rank 2.

\saa Let $ \caa $ be a two-sided cell of an affine Weyl group $ W $ of type $ \tilde{B}_2 $ and $ \tilde{G} _2$. Observation from Guilhot's partitions on two-sided cells shows that $ \caa $ can be written as the form\begin{equation}\label{eq:form}
\caa=\bigsqcup_{d\in D} B_d dU_d
\end{equation}
where $ B_d, U_d $ are subsets of $ W $, indexed by the elements $ d\in D $, each of  which is  an involution  in some finite parabolic subgroup. And the right cells in $ \caa $ are all of the form $ \Phi_{b,d}=bdU_d $ with $ b\in B_d $, $ d\in D $.

Let $ \mc{H}_{<\caa} $ be the two-sided ideal of the Hecke algebra $ \mc{H} $ spanned by the Kazhdan-Lusztig basis $ C_z $,  $ z<\caa $, where ``$ z<\caa $" means that $ z $ is in some two-sided cell $ \caa' $ with $ \caa'<\caa $ in the natural order on the set of two-sided cells. Then the main result of this paper can be stated as follows\begin{equation}\label{eq:in-dec}
C_{bdu}=E_bC_dF_u \mod{\mc{H}_{<\caa} } 
\end{equation}
 for all $ d\in D, b\in B_d,u\in U_d $, where $ E_b $ and $ F_u $ are determined by  equations $ C_{bd}=E_bC_d \mod{\mc{H}_{<\caa}}$, $ C_{du}=C_d F_u\mod{\mc{H}_{<\caa}}$. 
 
 As a corollary, the $ \mc{H} $-bimodule $ \mc{H}_{\leq\caa} /\mc{H}_{<\caa}$ (in the obvious meaning) is generated by $ C_d ,d\in D$. But a more important feature of this decomposition formula is that $ E_b $ and $ F_u $ are independent of each other. We will see that this feature is crucial for its application to conjectures P1-P15. 
 
 Assuming P4 and a weak form P8' of P8, and assuming $ \mb{a}(d)=\deg h_{d,d,d} $ for $ d\in D $, we can deduce  P1-P15 for $ \tilde{B}_2 $ and $ \tilde{G} _2$. The main tool is the decomposition formula. And some easy calculations are needed. If the $ \mb{a} $-value of each element of affine Weyl groups of type $ \tilde{B}_2 $ and $ \tilde{G} _2$ is known, then we can know immediately the validity of P4 and of the assumption $ \mb{a}(d)=\deg h_{d,d,d} $.  Thus computing $ \mb{a} $-functions of $ \tilde{B}_2 $ and $ \tilde{G} _2$ (with unequal parameters) is interesting and basic. But the we now have no method for it.

Note that  the form \eqref{eq:form} in general does not hold for any Coxeter group. One can check that the  Weyl group of type $ B_3 $ with equal parameter has given a counterexample. So it is a little surprise that \eqref{eq:form} holds for all the two-sided cells of affine Weyl groups of rank 2.

The results in this paper can be proved  for affine Weyl groups of type $ A_2 $  easily  in a similar way. So actually the affine Weyl group of rank 2 in this article  are used to refer to affine Weyl groups of type $ \tilde{B}_2 $ and $ \tilde{G}_2 $. 
 
\saa The organization of this paper is as follows. In Section 2, we recall some basic notions about Hecke algebras and fix some notations.  In Section 3, we recall Guilhot's partition of affine Weyl groups of rank 2 into two-sided cells and summarize some observations for latter use. In Section 4, we prove the decomposition formula under some assumptions, see Theorem \ref{th:dec}. In Section 5, we prove these assumptions are satisfied by affine Weyl groups of rank 2, and hence the decomposition formula holds in these cases. In Section 6, we apply the decomposition formula to conjectures P1-P15 (see Theorem \ref{th:ind}, \ref{th:commute}), and then reduce P1-P15 of affine Weyl groups of rank 2 to  P8' and computing Lusztig's $ \mb{a} $-function (see Theorem \ref{th:conjrk2}). At last, Appendix A and B contain some  computations in Section 5.

\section{Preliminaries}
\saa Let $(W,S)$ be a Coxeter group, where $S$ is  a finite set of generators of $W$ of order 2 satisfying braid relations. Let $\Gamma$ be an abelian group (written additively) equipped with an total order $ < $ such that $ \gamma+\lambda_1<\gamma+\gamma_2 $ for any $ \gamma \in \Gamma $ whenever $ \gamma_1<\gamma_1 $. Let $L:W\ra\Gamma$ be a positive weight function on $W$, i.e. $L$ is a map from $W$ to $\Gamma$ such that $L(ww')=L(w)+L(w')$ if $l(ww')=l(w)+l(w')$ and $L(w)>0$ if  $ w$ is not the neutral element $e $, where $l$ is the length function of Coxeter group $ W $.
Let $\mathcal{A}$ be the group algebra of $ \Gamma $. Then $\mc{A}=\mathbb{Z}[\Gamma]=\mathbb{Z}\{q^\gamma\,|\,\gamma\in\Gamma\}$ where $q$ is viewed as a symbol. We write $\Gamma_{<0}=\{\gamma\in\Gamma\,|\,\gamma<0 \}$ and $\mathcal{A}_{<0}=\mathbb{Z}\{q^\gamma\,|\,\gamma<0\}$; similarly for the notations $\Gamma_{>0}$, $\mathcal{A}_{>0}$, $\Gamma_{\leq0}$, $\mathcal{A}_{\leq0}$,etc..
Define the degree map $\deg:\mathcal{A}\ra\Gamma$ on $\mathcal{A}$ as  usual:\[
\deg(\sum_{\gamma\in\Gamma}a_\gamma q^{\gamma})=\max\{\gamma\,|\,a_\gamma\neq0 \}.
\]
This is well-defined since $\Gamma$ is totally ordered.

\saa\label{hecke} The Hecke algebra $\mathcal{H}$ associated with $(W,S,L)$ is an $\mathcal{A}$-algebra with $\mathcal{A}$-basis $\{T_w\,|\,w\in W \}$ and relations\[
T_wT_{w'}=T_{ww'}\text{ if } l(ww')=l(w)+l(w'),\]
\[(T_s+q^{-L(s)})(T_s-q^{L(s)})=0 \text{ if } s\in S.\]

There is a $\mathbb{Z}$-algebra involution $\bar{~}$ of $\mathcal{H}$, called bar involution,  such that $\bar{T}_w=T_{w^{-1}}^{-1}$, $\bar{q^{\gamma}}=q^{-\gamma}$ for  $w\in W$, $\gamma\in \Gamma$.  This bar involution is used to define the well known Kazhdan-Lusztig basis: there is a unique $\mathcal{A}$-basis of $\mathcal{H}$ such that
\[\bar{C}_w=C_w, \]
\[C_w\equiv T_w \mod{\mathcal{H}_{<0}},\] 
where $\mathcal{H}_{<0}:=\bigoplus_{w\in W}\mathcal{A}_{<0}T_w$.

Using Kazhdan-Lusztig basis, one can define partial order $ <_{\mc{L}} $ on $ W $, which is generated by the relation $ x\overset{\mc{L}}{\leftarrow} y$, where  $ x\overset{\mc{L}}{\leftarrow} y$ if there exist $ z\in W $ such that  $ C_x $ appears in the product $ C_zC_y $. The partial order $ <_{\mc{L}} $ induces naturally an equivalence relation   on   $ W $: $ x\overset{\mc{L}}{\sim} y$ if and only if $ x<_{\mc{L}} y$ and $ y<_{\mc{L}}x $. The corresponding equivalence class are called left cells. 

Similarly, one can define preoder $ <_\mc{R} $ on $ W $. Equivalently, $ x<_R y $ if and only if $ x^{-1}<_\mc{L}y^{-1} $. The corresponding equivalence relation is denoted by $ \overset{\mc{R}}{\sim} $ and the corresponding equivalence class are called right cells. 

At last, the preoder $ <_{\mc{LR}} $ is generated by $ <_\mc{L} $ and $ <_\mc{R} $.  The corresponding equivalence relation is denoted by $ \overset{\mc{LR}}{\sim} $ and the corresponding equivalence classes are called two-sided cells. A two-sided cell is usually denoted by $ \caa $ in this paper.

\saa We now introduce some notations for latter use.\begin{notation}\label{no:equal}
\begin{itemize}
\item [(i)] Obviously, the preorder $ \leq_{\mc{LR}} $ on $ W $ induces an partial order on the set of two-sided cells. We will denote this order by $ \leq $.
\item[(ii)] For $ x\in W $ and a two-sided cell $ \caa $, we  write $ x<\caa $ if $ x\notin\caa $ and $ x<_{\mc{LR}}w $ for any $ w\in \caa $. And we write $ x\leq\caa $ if $ x<_{\mc{LR}}w $ for any $ w\in \caa $.
\item [(iii)] Let $ \mc{H}_{<\caa}=\bigoplus_{z<\caa}\mc{A} C_z$. Then $  \mc{H}_{<\caa} $ is a two-sided ideal of $ \mc{H} $. We have similar meaning for $ \mc{H}_{\leq \caa} $.

\item [(iii)] For $ a,b\in\mc{H} $,
\begin{itemize}
\item we write $ a\equiv b $ for  ``$a-b\in{\mc{H}_{<0}}$", or equivalently  for ``$ a=b\mod{\mc{H}_{<0}} $'';
\item We write  $a\ed b$  for  ``$a-b\in{\mc{H}_{<\caa}}$", or equivalently   ``$ a=b\mod{\mc{H}_{<\caa}} $",  or ``$ a=b $ in $ \Hm/\Hm_{<\caa} $";
 \item  We write $a\edq b$ for  ``$a-b\in \mc{H}_{<0}+\mc{H}_{<\caa}$",  or equivalently  ``$ a=b\mod{\mc{H}_{<0}+\mc{H}_{<\caa}} $", where $  \mc{H}_{<0}+\mc{H}_{<\caa}=\{a+b\mid a\in \Hm_{<0}, b\in \Hm_{<\caa} \} $.
\end{itemize}
\item[(iv)] For $ \gamma\in\Gamma $, we write $ \xi_\gamma=q^\gamma-q^{-\gamma},  \eta_\gamma=q^\gamma+q^{-\gamma}  $.
\end{itemize}
\end{notation}

The following lemma is useful for the proof and calculation. 
\begin{lemma}\label{lem:equal}
If $ a,b\in\mc{H} $  are both bar invariant in $ \Hm/\Hm_{<\caa} $ and  $ a= b \mod{\mc{H}_{<0}+\mc{H}_{<\caa}}$, then $ a=b\mod{\mc{H}_{<\caa}} $. Using notations in  \ref{no:equal}, if $ a\ed\bar{a} $, $  b\ed\bar{b} $ and $ a\edq b$ then  $ a\ed b  $.
\end{lemma}
\textit{Proof.} 
The lemma is equivalent to that if $a-b \in \mc{H}_{<0}+\mc{H}_{<\caa}$ and is bar invariant in  $ \Hm/\Hm_{<\caa} $, then $ a-b\in \mc{H}_{<\caa} $. Hence we only need to prove the following claim \begin{equation}\label{eq:cm1}
\text{If  }a\in\mc{H}_{<0}+\mc{H}_{<\caa} \text{ and  } \bar{a}-a\in\mc{H}_{<\caa}, \text{ then } a\in \mc{H}_{<\caa}. 
\end{equation}

By linear independence, we have unique  $ \alpha_y,\alpha_z \in\mc{A}$ such that \[
a=\sum_{y\nless \caa}\alpha_yC_y+\sum_{z< \caa}\alpha_zC_z. 
\]Let $a'=\sum_{y\nless \caa}\alpha_yC_y$. Since $a\in\mc{H}_{<0}+\mc{H}_{<\caa}$,  we have $\alpha_y\in\mc{A}_{<0}$ and hence \begin{equation}\label{eq:less}
 a'\in \mc{H}_{<0}.
\end{equation} 
On the other hand,  since $ \bar{a}-a \in\mc{H}_{<\caa}$, we have $ \bar{a'}-a'\in\mc{H}_{<\caa}\cap(\bigoplus_{y\nless \caa}\mc{A}C_y)=0 $, and hence  \begin{equation}\label{eq:barrr}
\bar{a'}=a'.
\end{equation}
Combining \eqref{eq:less} and  \eqref{eq:barrr}, we get $a'=0$, see \cite[\S 5(e)]{lusztig2003hecke}.  This completes the proof of the claim  \eqref{eq:cm1}. \qed

\saa \label{sec:map}There is an anti-involution $ {~}^\flat $ of $ \mc{A} $-algebras on the Hecke algebra $ \mc{H} $ such that $ T_w^{\flat}=T_w^{-1} $. It is obvious that $ (hC_w)^\flat=C_{w^{-1}}h^\flat $ for any $ h\in\mc{H} $.

Define the $ \mc{A} $-linear map $ \tau:\mc{H}\ra\mc{A} $ such that $ \tau(T_x)=\delta_{x,e} $. It is well known that $ \tau(T_xT_y)=\delta_{x,y^{-1}} $ and $ \tau(hh')=\tau(h'h) $ for any $ h,h'\in\mc{H} $. One can check that $ \tau(C_xC_y)=\delta_{x,y^{-1}}\mod{\mc{H}_{<0} }$.

\saa Let $P_{y,w}$, $y\in W$ be the elements in $\mathcal{A}_{\leq0}$ such that $C_w=\sum_{y\in W}P_{y,w}T_y$. Define $\Delta(z)$ to be an element of $\Gamma_{\geq0}$ and $n_z\in\mathbb{Z}$ such that $P_{e,z}$ ($e$ is the neutral element of $W$) is of the form\[
P_{e,z}=n_zq^{-\Delta(z)}+\text{lower degree terms}.
\] 

Let $h_{x,y,z}\in \mathcal{A}$ be defined by $C_{x}C_{y}=\sum_{z}h_{x,y,z}C_z$. Then   the $\textbf{a}$-function $\textbf{a}:W\ra\Gamma\cup\{\infty \}$ on $W$ is defined by\[
\textbf{a}(z):=\sup\{\deg(h_{x,y,z})\,|\,x,y\in W\}.
\]
If $\textbf{a}(z)$ is finite, $h_{x,y,z}$ can be written as the form
\[
h_{x,y,z}=\gamma_{x,y,z^{-1}}q^{\textbf{a}(z)}+\text{ lower degree terms}, ~\gamma_{x,y,z^{-1}}\in\mathbb{Z}
\] Recall that   the $\textbf{a}$-function is always bounded for  affine Weyl groups (with any parameters), see \cite{lusztig1985cellsI}.

Define $\mathcal{D}:=\{z\in W\,|\,\textbf{a}(z)=\Delta(z)\}.$

Now we recall Lusztig's conjectures P1-P15 (\cite[\S14.2]{lusztig2003hecke}).  The following is actually  a reformulation of  P1-P15.

\saa\label{P1-P15} \textbf{Conjectures} $ (P1)_{\leq\caa}$-$(P15)_{\leq\caa} $. Let $ (W,S,L) $ be a Coxeter group with positive weight function  and with bounded $\textbf{a}$-function. Fix a two-sided cell $ \caa $. Then

$ (P1)_{\leq\caa} $. For any  $z\leq\caa$, we have $\mb{a}(z)\leq\Delta(z)$.

$ (P2)_{\leq\caa} $.  If  $q\in\mathcal{D}_{\leq\caa}=\mathcal{D}\cap\{x\mid x\leq\caa\}$, and   $x,y\leq\caa$ are such that $\gamma_{x,y,q}\neq0$, then  $x=y^{-1}$.

$ (P3)_{\leq\caa} $. For any  $y\leq\caa$, there exists  uniquely  $q\in\mathcal{D}_{\leq \caa}$ such that  $\gamma_{y^{-1},y,q}\neq0$.

$ (P4)_{\leq\caa} $. If  $ z,z'\leq \caa $ and  $z'\leq_{ {LR}}z$, then $\mb{a}(z')\geq\mb{a}(z)$. In particular if $ z'\sim_{LR} z $ then $ \mb{a}(z)=\mb{a}(z') $.

$ (P5)_{\leq\caa} $.  If $q\in\mathcal{D}_{\leq\caa}$, $ y\leq \caa $, $\gamma_{y^{-1},y,q}\neq0$, then $\gamma_{y^{-1},y,q}=n_q=\pm1$.

$ (P6)_{\leq\caa} $. For any $q\in\mathcal{D}_{\leq\caa}$, we have  $q^2=1$.

$ (P7)_{\leq\caa} $. For any  $x,y,z\leq\caa$, we have $\gamma_{x,y,z}=\gamma_{y,z,x}$.

$ (P8)_{\leq\caa} $. For any $x,y,z\leq\caa$, $\gamma_{x,y,z}\neq0$ implies that  $x\sim_{\mathcal{L}}y^{-1}$, $y\sim_{ \mc{L}}z^{-1}$, $z\sim_{ \mc{L}}x^{-1}$.

$ (P9)_{\leq\caa} $.  For any  $z,z'\leq\caa$, if  $z'\leq_{\mc {L}}z$, $\mb{a}(z')=\mb{a}(z)$, then $z'\sim_{ \mc{L}}z$.

$ (P10)_{\leq\caa} $.  For any  $z,z'\leq\caa$, if  $z'\leq_{\mc{R}}z$, $\mb{a}(z')=\mb{a}(z)$, then $z'\sim_{ \mc{R}}z$.

$ (P11)_{\leq\caa} $.  For any $z,z'\leq\caa$, if  $z'\leq_{ \mc {LR}}z$, $\mb{a}(z')=\mb{a}(z)$, then $z'\sim_{ \mc{LR}}z$.

$ (P12)_{\leq\caa} $. Let $I\subset S$. If  $y\in W_I$ and $ y\leq \caa $, then the $ \mb{a} $-value of $ y $ in  $W_I$ is the same as that in $W$. 

$ (P13)_{\leq\caa} $. Any left cell  $\Theta\subset \{x\mid x\leq\caa \}$ contains  a unique element $ q $ in  $\mathcal{D}_{\leq \caa}$. And for each $y\in\Theta$ we have $\gamma_{y^{-1},y,q}\neq0$. Similar property holds for right cells.

$ (P14)_{\leq\caa} $. For any $z\leq \caa$,  we have  $z\sim_{ \mc{LR}}z^{-1}$.

$ (P15)_{\leq\caa} $. For any  $w,w'\in W$ and any $ x,y $ such that  $ x\sim_{\mc{LR}} y \leq\caa$, we have\[
\sum_{z\in W} h_{x,w',z}\otimes h_{w,z,y}=\sum_{z\in W}h_{z,w',y}\otimes h_{w,x,z} \in  \mc{A}\otimes_\mathbb{Z}\mc{A}.
\]

 We will need a weak version of $ (P8)_{\leq\caa} $:

$ (P8')_{\leq\caa} $. For any $x,y,z\leq\caa$, $\gamma_{x,y,z}\neq0$ implies that $y\sim_{ \mc{L}}z^{-1}$.

If $ \caa$ is the highest two-sided cell $\{e\} $,  we simply  denote $ (P1)_{\leq\caa}$-$ (P15)_{\leq\caa}$   by P1-P15, respectively.

\section{Partition into cells}\label{p-c}

In this section we will summarize the partition of  two-sided cells of affine Weyl group of type $ \tilde{B}_2 $, and  $ \tilde{G}_2 $ which is due to Guilhot, see \cite{guilhot2010rank2,guilhot0810.5165}. Then we conclude some easy properties  of these cells for latter use.

\begin{figure}
\centering
\includegraphics[width=0.3\linewidth]{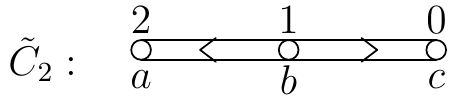}
\caption{Dynkin diagram of $\tilde{C}_2$.}
\label{fig:ad-eps-converted-to}
\end{figure}

\saa
Assume first that $W$ is an affine Weyl group of type $\tilde{C}_2$ with the set of simple reflections $S=\{s_0,s_1,s_2 \}$ and the relations $s_0s_1s_0s_1=s_1s_0s_1s_0$, $s_2s_1s_2s_1=s_1s_2s_1s_2$ and $s_0s_2=s_2s_0$. The weight function $L:W\ra\Gamma$ is given by $L(s_2)=a$, $L(s_1)=b$, $L(s_0)=c$ with $a,b,c$ all in $\Gamma_{>0}$. See Figure \ref{fig:ad-eps-converted-to}.

We often use a sequence of 0,1,2 to denote an element of $W$. For example, if $w=s_0s_1s_0s_1s_2$, we abbreviates $w$ by $01012$, $T_w$ by $T_{01012}$, and $C_w$ by $C_{01012}$. The empty sequence denotes the neutral element $e$ of $ W $. 

If $A,B\subset W$, we write $AB$ for the set $\{xy\,|\,x\in A,y\in B \}$; similar for $ABC$ with $A,B,C\subset W$. If we write $AB$, we automatically have $l(xy)=l(x)+l(y)$ for any $x\in A, y\in B$ unless otherwise specified.

We denote by $\leq_D$ the Duflo oder on $ W $. In other words, \[x\leq_D y\text{ if and only if }l(x^{-1}y)=l(y)-l(x).\] If $p\in W$ , we write $$ U(p)=\{ w\in W\,|\, w\leq_D p^k \text{ for some }k\in\mathbb{N} \}. $$

\saa\label{pc2} We classify the two-sided cells of the affine Weyl group of type $\tilde{C}_2$ for various positive parameters ($a,b,c>0$ and $a\geq c$) into the following types(\cite{guilhot0810.5165}).  

\begin{itemize}
\item [(1)]Let $\caa= BdU$, $ U=U(p) $ with
  $B$, $d$, $p$ taking elements in the following table.
  
\begin{tabular}{|c|c|c|c|c|}
\hline cases & parameter condition & $B$ &  $d$& $p$ \\ 
\hline i & $a-c>2b$ & $\{e,1,01,101\}$ & 02 & 1012  \\ 
\hline ii&$0<a-c<2b$  & $\{e,2,12,,012\}$ & 1010 & 2101  \\ 
\hline iii& $|a-c|<b,a+c>b$ &  $\{e,1,01,,21\}$ & 02 &  102\\ 
\hline iv & $a+c<b$ &  $\{e,0,2,02\}$&1  &  021\\ 
\hline v &  $a-c>b$ & $\{e,0,10,010\}$& 212  & 012  \\ 
\hline vi& $a>c,a+c>2b$ & $\{e,1,01,101\}$ & 2 &  1012\\ 
\hline vii& $a>c,a+c<2b$  &  $\{e,2,12,,012\}$& 101&  2101\\ 
\hline
\end{tabular}

Then $ \caa $ is a two-sided cell. Set $ D=\{d\} $.
\item [(2)] Let $\caa=BdB^{-1}$ with
 $B$, $d$ taking elements in the following tables.

\begin{tabular}{|c|c|c|c|}
\hline cases & condition & $B$ &  $d$\\ 
\hline i& $a<b, c<b,a+c>b$ &$\{e,0,2\}$ &1 \\
\hline  ii &$c<a<b$ &$\{e,0\}$& {121}\\
\hline  iii &$a>b>c$ &$\{e,0\}$& 1\\
\hline  iv &$b<a-c<2b$ &$\{e,1,01\}$& 02\\
\hline  v &$a>b,a+c<2b$ &$\{e,1,01\}$& 2\\
\hline  vi &$a>c>b$ &$\{e,1\}$& 0\\
\hline
\end{tabular}

 Then $ \caa $ is a two-sided cell. Set $ D=\{d\} $ and  $U=B^{-1}$.

\item [(3)]   Let $\caa=B_1d_1U_1\cup B_2d_2U_2$,
 $U_1=U(p_1)$, $U_2=U(p_2)$ with  $B_1$, $d_1$, $p_1$ and  $B_2$, $d_2$, $p_2$ taking elements in the following table.

\begin{tabular}{|c|c|c|c|c|c|c|c|}
\hline cases & conditions  & $B_1$ & $d_1$& $p_1$ & $B_2$ & $d_2$ & $p_2$  \\ 
\hline i& $a-c=2b$ & $\{e,1,01\}$ & 02 & 1012 & $\{e\}$ & 1010  & 2101  \\ 
\hline ii & $a+c=b$ &$\{e,0,2\} $ & 1 &021  & $\{e\} $& 02  & 102  \\ 
\hline iii & $a-c=b$ &$\{e,1,01\}$  & 02 &  102 & $\{e\}$ & 212 & 012 \\ 
\hline  iv& $a>c,a+c=2b$ &$\{e,1,01\}$  & 2  & 1012 &$\{e\} $ & 101  & 2101  \\ 
\hline v & $a=c,a>b$  &$\{e,1\}$  & 2  & 1012  & $\{e,1\}$  & 0  & 1210 \\ 
\hline  vi & $a=c,a<b$ & $\{e,2\}$ &  101 &  2101 & $\{e,0\}$  & 121  & 0121  \\ 
\hline 
\end{tabular} 

Then $ \caa $ is a two-sided cell. Set $ D=\{d_1,d_2\} .$
\item [(4)] Exotic cases: 
\begin{itemize}
\item [(i)]If $a=b=c$,  set $B_1=\{e\}$, $d_1=1$, $U_1=U(0121)\cup U(2101)$, set $B_2=\{e\}$, $d_2=2$, $U_2=U(1012)\cup\{12\}$, and set $B_3=\{e\}$, $d_3=0$, $U_3=U(1210)\cup\{10\}$. Then $\caa=B_1d_1U_1\cup B_2d_2U_2\cup B_3d_3U_3$ is a  two-sided cell. Set $D=\{0,1,2\}$.
\item [(ii)]If $a>b=c$, then $$\caa=\{1,10,101\}\cup\{0,01,010\}$$ is a  two-sided cell. Set $B_1=\{e\}$, $d_1=1$, $U_1=\{e,0,01\}$ and set $B_2=\{e\}$, $d_2=0$, $U_2=\{e,1,10\}$. Then $\caa=B_1d_1U_1\cup B_2d_2U_2$. Set $D=\{d_1,d_2\}$.
\item [(iii)]If $a=b>c$, then $$\caa=\{1,10,12,121,1210\}\cup\{01,010,012,0121,01210\}\cup\{2,21,210,212\}$$ is  a two-sided cell. Set $B_1=\{e,0\}$, $d_1=1$, $U_1=\{0,2,21,210\}$, and set $B_2=\{e\}$, $d_2=2$, $U_2=\{e,1,10,12\}$. Then $\caa=B_1d_1U_1\cup B_2d_2U_2$. Set $D=\{d_1,d_2\}$.
\end{itemize}
\item [(5)]All the cells with one element are listed as follows:
\begin{itemize}
\item $ \{1\} $ when $ a,c>b $;
\item $ \{010\} $ when $ a,c>b $;
\item $ \{1010\} $ when $ a-c>2b $;
\item $ \{212\} $ when $ b,c<a $ and $ a<b+c $;
\item $ \{0\} $ when $ c<b $ and $ a\geq c $;
\item $ \{101\} $ when $ a+c>2b $ and $ c<b $;
\item $ \{2\} $ when $ a\geq c $ and $ a<b $;
\item $ \{02\} $ when $ a+c<b $.
\end{itemize}
\item [(6)] The lowest two-sided cell $ \caa $:
\begin{itemize}
\item if $ a>c $, then  \[ \caa=BdU \] with $ B=\{b\mid b^{-1}\leq_D 010210 \} $, $d=1212$, and $ U=\{ u\in W\mid l(1212u)= l(u)+l(1212)\} $;
\item if $ a=c $, then \[ \caa=B_1d_1U_1\cup B_2d_2U_2 \] with $ d_1=1212 $, $ d_2=1010 $, $ U_1=\{ e,0,10,210 \} $, $ B_2=\{e,2,12,012 \} $, $ U_1=\{ u\in W\mid l(1212u)= l(u)+l(1212)\} $, and $ U_2=\{ u\in W\mid l(1010u)= l(u)+l(1010)\} $.
\end{itemize}
\end{itemize}

These are all   possible two-sided cells of $ \tilde{C}_2 $.

\begin{figure}
\centering
\includegraphics[width=0.3\linewidth]{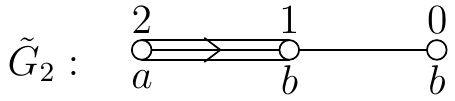}
\caption{Dynkin diagram of $\tilde{G}_2$.}
\label{fig:adg-eps-converted-to}
\end{figure}

\saa\label{g2} Assume now that $W$ is the affine Weyl group of type $\tilde{G}_2$ with   generators $ s_0,s_1,s_2 $ and  relations $s_1s_2s_1s_2s_1s_2=s_2s_1s_2s_1s_2s_1$, $s_1s_0s_1=s_0s_1s_0$, and $s_0s_2=s_2s_0$. The weights are $L(s_2)=a$, $L(s_1)=b=L(s_0)$ with $a,b\in\Gamma_{>0}$. See Figure \ref{fig:adg-eps-converted-to}.

We also abbreviate an element of $ W $ using the sequences of 0,1,2  associated with its reduced expressions. 
We now describe   the two-sided  cells of $\tilde{G}_2$.
\begin{itemize}
\item [(1)]
Let $\caa=BdU$, $U=U(p)$ with   $B$, $d$, $p$ taking elements in the following table.

\begin{tabular}{|c|c|c|c|c|}
\hline cases & condition & $B$ &  $d$& $p$ \\ 
\hline i & $2a>3b$ & $\{e,0,10,210,1210,01210\}$ & 21212 & 01212  \\ 
\hline ii & $2a<3b$ & $\{e,1,21,121,0121,2121\}$ & 02 & 12102\\
\hline iii & $a>2b$ & $\{e,1,01,21,121,0121\}$ & 02 & 102\\
\hline iv &$a<2b$ &$\{e,2,12,212,1212,01212\}$ & 101 & 210\\
\hline v &$2a>3b,a<2b$& $\{e,1,21,121,0121\}$ & 02 & e\\
\hline vi &$a<b$& $\{e,0\}$& 12121& e\\
\hline
 \end{tabular} 
 
 Then $ \caa $ is a two-sided cell. Set $ D=\{d\} $.
 \item [(2)]  Let $\caa=B_1d_1U_1\cup B_2d_2U_2$, $U_1=U(p_1)$, $U_2=U(p_2)$ with $B_1$, $d_1$, $p_1$ and  $B_2$, $d_2$, $p_2$ taking values in the following table.
 
 \begin{tabular}{|c|c|c|c|c|c|c|c|}
 \hline cases & condition  & $B_1$ & $d_1$& $p_1$ & $B_2$ & $d_2$ & $p_2$  \\ 
 \hline i & $2a=3b$ & $\{e,1,21,121,0121\}$ & 02 & 12102 & $\{e\}$ & 21212  & 01212  \\ 
 \hline ii & $a=2b$ &$\{e,1,21,121,0121\} $ & 02 & 102  & $\{e\} $& 101  & 210  \\ 
 \hline 
 \end{tabular}  Then $ \caa $ is a two-sided cell. Set $D=\{d_1,d_2\}$.
 \item [(3)] Exotic cases.
 \begin{itemize}
 \item [(i)] If $a>b$, then  $$\caa=\{1,10\}\cup\{0,01\}$$ is a two-sided cell; $D=\{1,0\}$. Set $B_1=B_2=\{e\}$, $d_1=1$, $d_2=0$, $U_1=\{e,0\}$, $U_2=\{e,1\}$. Then $\caa=B_1d_1U_1\cup B_2d_2U_2$.
 \item [(ii)]  If $a<b$, then  \begin{multline*}
 \caa=\{21,210,212,2121,21210,21212\}\cup\{1,12,121,1212,10,1210\}\\\cup\{0,01,012,0121,01210,01212\}
 \end{multline*} is a two-sided cell; $D=\{1,0\}$. Set $B_1=\{e,2\}$, $d_1=1$, $U_1=\{e,2,0,21,212,210\}$ and set $B_2=\{e\}$, $d_2=0$, $U_2=\{e,1,12,121,1210,1212\}$. Then $\caa=B_1d_1U_1\cup B_2d_2U_2$.
 \item [(iii)] If $a=b$, then \begin{multline*}
 \caa=\{1,10,12,121,1210,1212,12121,121210\}\\\cup\{0,01,012,0121,01210,01212,012121,0121210\}\\\cup\{2,21,210,212,2121,21210,21212\}
 \end{multline*} is a two-sided cell; $D=\{1,0,2\}$.
 Then set $d_1=1$, $d_2=0$, $d_3=2$. Then $\caa=B_1d_1U_1\cup B_2d_2U_2\cup B_3d_3U_3$ where $B_1=B_2=B_3=\{e\}$.
 \item [(iv)]  If $a>b$, then
 \[
 \caa=\{e,1,01\}\{2\}\{e,12\}\{e,1,10 \}.
 \]
 is a two-sided cell. Set $ d=2 $, $ B=\{e,1,01\} $, $ U=\{e,1,10,12,121,1210\} $, $ D=\{2\} $. Then $\caa=BdU  $.
 \end{itemize}
 \item [(4)] The two-sided  cells with one element:
 \begin{itemize}
 \item $ \{2\} $ when $ a<b $;
 \item $ \{21212\} $ when $ 2b<2a<3b $;
 \item $ \{101\} $ when $ a>2b $.
 \end{itemize}
 \item [(5)] The lowest two-sided cell: 
 \[\caa=BdU\]
 where $ B=\{b\mid b^{-1}\leq_D 0121201210  \} $, $ d=121212 $,  	$ U=\{u\in W\mid l(121212u)=l(121212)+l(u) \} $.
\end{itemize}

These are all   possible two-sided cells of $ \tilde{G}_2 $.

\saa We conclude this section with following observations.
\begin{lemma}\label{form}
Let $\caa$ be a two-sided cell of an affine Weyl group $W$ of type $\tilde{C}_2$ or $\tilde{G}_2$ with positive weight function. Then 
\begin{itemize}
\item [(i)]$\caa$ can be written in the form $$ \caa=\bigcup_{d\in D}B_ddU_d $$ where $D$, $B_d$, $U_d$ ($d\in D$) are subsets  of $W$ which have been listed in \ref{pc2}, \ref{g2} case by case. And this form satisfies\begin{itemize}
\item For each $ d\in D $, $ b\in B_d $ and  $ u\in U_d $, we have $ l(bdu)=l(b)+l(d)+l(u) $.
\item  $ e \in B_d,e\in U_d$, and $ B_d^{-1}\subset U_d $.
\item Each  $ d\in D $ is  an involution element of some finite parabolic subgroups of $ W $ such that $ l(su)=l(u)+1 $ for any  $ u\in U_d $, and $ s\in S $ with $s\leq d $.
\end{itemize}

\item [(ii)]If  $ w\notin U_d$ satisfying  $ l(sw)=l(w)+1 $  for $ s\in S $ with  $ s\leq d  $,   then $ dw\notin \caa $.

\item [(iii)]The right cells in $\caa$ are all of the form \[
\{bdu\,|\,u\in U_d \}, d\in D,b\in B_d.
\]
\end{itemize}
\end{lemma}\textit{Proof.} (i) is just a direct observation from subsection \ref{pc2}, \ref{g2}. (ii) follows from \cite[Rem,6.7]{guilhot2010rank2}.  (iii) is one of the conclusion in \cite{guilhot2010rank2}.\qed

\section{Decomposition formula}

In this section, we need the following assumptions on a two-sided cell $ \mb{c} $ of a Coxeter group $ W $.

\begin{assumption}\label{ass}
Let $ \mb{c} $ be a two-sided cell of $ W $. We assume that there exist subsets $ D $, $ B_d $, $ U_d $ $ (d\in D) $ of $ \mb{c} $ such that
\begin{itemize}
\item [(i)]\begin{itemize}
\item [(a)] $\caa^{-1}= \caa=\bigsqcup_{d\in D}B_ddU_d $.
 \item  [(b)]For each $ d\in D $, $ b\in B_d $ and  $ u\in U_d $, we have $ l(bdu)=l(b)+l(d)+l(u) $.
 \item [(c)]  The neutral element $ e \in B_d,U_d$, and $ B_d^{-1}\subset U_d $.
\item [(d)] If  $ w\notin U_d$ satisfying  $ l(sw)=l(w)+1 $  for any simple reflection $ s\leq d  $,   then $ dw\notin \caa $.
\end{itemize}

\item [(ii)] Each  $ d\in D $ is  an involution element of some finite parabolic subgroups of $ W $ with the following properties:
\begin{itemize}
\item [(a)]$ l(su)=l(u)+1 $ for any $ s\leq d $, $ u\in U_d $.
\item [(b)] For any  $ s\leq d  $, we have  $ T_sC_d\in \mc{A}C_d+\Hm_{<\caa} $.
\item [(c)] $ h_{d,d,d}\neq0 $.
\end{itemize}

\item [(iii)] For any $ d\in D $, the set $ dU_d $ is a right cell of $ W $.
\item [(iv)]  For any $ d\in D $, $ b\in B_d $, $ u\in U_d $,  we have  $$  \T{b}C_d\T{u}=\T{bdu} \mod{\mc{H}_{<0}+\mc{H}_{<\caa}}. $$ 
\end{itemize}\end{assumption}

It follows  from Lemma \ref{form} that  Assumption \ref{ass} has been hold for affine Weyl groups of rank 2, except Assumption \ref{ass} (ii,b-c) and (iv). We will verify them with some calculations in the next section.

The main result of this section can roughly be stated as follows. Under the assumption \ref{ass}, we have product decomposition:
$$ h_{d,d,d}C_{bdu}=C_{bd}C_{du}\mod{\mc{H}_{<\caa}} $$
for any $ d\in D, b\in B_d, u\in U_d $.

\begin{lemma}Keep  Assumption \ref{ass}. \label{lem:indepen}
\begin{itemize}
\item [(i)]For any $ d\in D $, $ u\in U_d $,  we have $ C_d\T{u}=\T{du}\mod{\Hm}_{<0} $.
\item [(ii)] For any  $ d\in D $, we have  $ C_dC_d=h_{d,d,d}C_d\mod{\Hm}_{<\caa} $.
\item [(iii)]The set $\{ C_dT_u\mid u\in U_d\}  $ is a $ \mc{A} $-linearly independent subset of    $ \mc{H}/\mc{H}_{<\caa} $.
\end{itemize}
\end{lemma}
\textit{Proof.} The statement (i) follows directly from Assumption \ref{ass} (ii,a). The statement (ii) follows directly from Assumption \ref{ass} (ii,b). To prove (iii), we only need to prove that \begin{equation*}
 \{C_dT_u\mid u\in U_d \}\cup\{C_z\mid z<\caa\} \text{ is liearly independent in }  \mc{H},
\end{equation*}
which follows from the fact that
 $$ C_dT_u \in T_{du}+\sum_{y<du}\mc{A}T_y,\text{ and }C_z\in T_z+\sum_{y<z}\mc{A}T_y .$$\qed

\begin{defn}
Fix an element $ d\in D $. For any $ w\in W $, we say $ y<_U w $ if and only if $ y<w $ in Bruhat order and $ y\in U_d $. And we write $ y\leq_U w $ if $ y<_Uw $ or $ y=w $. In other words, for any $ w\in W $, we have $\{y\mid y\leq_U w \}=\{y\in U_d\mid y< w\}\cup\{w\} $.
\end{defn}

\begin{lemma}\label{lem:loop}
Fix a two-sided cell $ \mb{c} $ as the Assumption \ref{ass} and an element $ d\in D $. Let $ w\in W $ satisfy the condition that for any  element $ s\in S $ such that $ s\leq d $ we have $ l(sw)=1+l(w) $ (the elements of  $ U_d $ satisfy this condition). Then \begin{itemize}
\item [(i)]There exists $ r_{y,w} \in\mc{A}$ such that \[ \overline{C_dT_w}=C_dT_w+\sum_{y<_U w}r_{y,w}C_dT_{y}\mod{\mc{H}_{<\caa}} .\]

\item [(ii)] There exists a unique element $ F_w\in\mc{H} $ such that\begin{itemize}
\item [(a)]$F_w=T_w+\sum_{y<_Uw}p_{y,w}T_{y}$ with $ p_{y,w}\in\mc{A}_{<0} $.
\item [(b)] $ C_dF_w=C_{dw} \mod{\mc{H}_{<\caa}}$.
\end{itemize}\item [(iii)]If  $w\notin U_d$, then $ C_dT_w\in \sum_{y<_Uw}\mc{A}C_dT_y+\mc{H}_{<\caa}$.
\end{itemize} 
\end{lemma}
\textit{Proof.} We use induction on the length $ n=l(w) $ of $ w $. The proof forms an interesting loop. The lemma for $ n=0 $ is obvious. We will use the notation $ (iii)_{<n} $ to express that  the statement (iii) holds for all  $ w $ with length $ <n $. We will prove that $ (iii)_{<n} $ implies $ (i)_{n} $; $ (i)_{\leq n} $ implies $ (ii)_{n} $; and $ (ii)_{n} $ implies $ (iii)_n $. Then we can conclude the lemma.

$ (iii)_{<n} \implies (i)_{n} $. It is well-known that $ \overline{T_w}=T_w +\sum_{y'<w}R_{y',w}T_{y'}$ with $ R_{y',w}\in\mc{A} $. By the Assumption \ref{ass} (ii,b), we have  \[\overline{C_d\T{w}}\in C_dT_w+\sum_{y}\mc{A}C_dT_y+\Hm_{<\caa},\] where $y$ takes over the elements such that \[y<w,\text{ and } l(sy)=1+l(y)\text{ for any simple relection } s\leq d. \] Applying $ (iii)_{<n-1} $ to $y$,  we immediately get $ (i)_n $.

$ (i)_{\leq n} \implies(ii)_{n} $. The basic idea follows from the construction of Kazhdan-Lusztig basis. 
We first prove that\begin{itemize}
\item [(*)] there exists a unique element $ F_u\in\mc{H} $ such that   (a)  holds and  $ C_dF_u $ is  invariant in  $ \mc{H}/\mc{H}_{<\caa} $ under the bar involution $ \bar{~} $.  
\end{itemize}

Let $ w' $ be such that $ l(w')\leq n $ and such that $ l(sw')=l(w')+1 $ for any $ s\leq d $ with $ s\in S $. By $ (i)_{\leq n} $, we have \[\overline{C_dT_{w'}}=\sum_{y\leq_U{w'}}r_{y,w'}C_dT_{y} \mod{\mc{H}/\mc{H}_{<\caa}} ,\] where we  convention that $r_{w',w'}=1$. Then\[
 C_dT_{w'}\ed\sum_{y\leq_U w'}\bar{r}_{y,w'}\overline{C_dT_{y}} \ed\sum_{y'\leq_U w'}\left(\sum_{y:y'\leq_U y\leq_U w'}{r}_{y',y}\bar{r}_{y,w'}\right) C_dT_{y'}.
\] Recall from Notation \ref{no:equal} that  $`` \ed" $ means ``equals by modulo $ \mc{H}_{<\caa} $". By linear independence (see Lemma \ref{lem:indepen}(iii)), we obtain  \begin{equation}\label{eq:sel}
  \sum_{y:y'\leq_U y\leq_U w'}{r}_{w',y}\bar{r}_{y,w'}=\delta_{y',w'}.
  \end{equation}
(Note that when $ y'=w' $, the equation   holds by our convention $ r_{w',w'} =1$.)

Write $ F_w=\sum_{y\leq_U w}p_{y,w}T_y $ and convention that $ p_{w,w}=1 $. We prove now claim (*) by solving equations. It is easy to see that \[
 \overline{C_{d}F_{w}}\ed \sum_{y'\leq_Uw}(\sum_{y:y'\leq_Uy\leq_Uw}r_{y',y}\bar{p}_{y,w})C_dT_{y'}.
\]
By  Lemma \ref{lem:indepen}(iii), we see that to require  $C_dF_u$ is invariant in $ \mc{H}/\mc{H}_{<\caa} $ under bar involution is equivalent to require that \begin{equation*}
 \sum_{y:y'\leq_Uy\leq_Uw}r_{y',y}\bar{p}_{y,w}=p_{y',w},\forall y'\leq_U w
  \end{equation*}
  and this is further equivalent to \begin{equation}\label{eq:sol}
p_{y',w}-\bar{p}_{y',w}=\sum_{y:y'<_Uy\leq_Uw}r_{y',y}\bar{p}_{y,w}.
  \end{equation} 
We  solve $ p_{y',w} $ by induction on $ l(w)-l(y') $. If $ l(w)-l(y')=0 $ then \eqref{eq:sol} agrees with our convention $ p_{w,w}=1 $.  Assume that we have known $ p_{y,w} $ for all $ y $ satisfying $ y'<_U y\leq_Uw $. Then $ p_{y,w}= \sum_{v:y\leq_Uv\leq_Uu}r_{y,v}\bar{p}_{v,w} $. So	\begin{align*}
   \overline{\sum_{y:y'<_Uy\leq_Uw}r_{y',y}\bar{p}_{y,w}}&=\sum_{y:y'<_Uy\leq_Uw}\bar{r}_{y',y}(\sum_{v:y\leq_Uv\leq_Uw}r_{y,v}\bar{p}_{v,w})\\
   &=\sum_{v:y'<_Uv\leq w}(\sum_{y:y'<_Uy\leq_Uv}\bar{r}_{y',y}r_{y,v})\bar{p}_{v,w}\\
   &=\sum_{v:y'<_Uv\leq w}(-r_{y',v})\bar{p}_{v,w}\text{ (by \eqref{eq:sel}) }\\
   &=-\sum_{y:y'<_Uy\leq_Uw}r_{y',y}\bar{p}_{y,w}
   \end{align*}
In other words, the right hand side of \eqref{eq:sol} is anti-invariant under bar involution. Then we can see that $ p_{y,w} $ is just the  part of the negative degrees of the Laurent polynomial $ \sum_{y:y'<_Uy\leq_Uw}r_{y',y}\bar{p}_{y,w} $.  So $ p_{y',w} $ uniquely determined by \eqref{eq:sol} and hence (*) holds.

By Lemma \ref{lem:indepen}(i) and the assumption on $ w $, we have \[
C_dF_w=T_{dw}=C_{dw}\mod{\mc{H}_{<0}}.
\]Since $ C_dF_w $ is  invariant in  $ \mc{H}/\mc{H}_{<\caa} $ under the bar involution, we have $ C_dF_w =C_{dw}\mod{\mc{H}_{<\caa}}$, see Lemma \ref{lem:equal}. This proves $ (ii)_n $.

$ (ii)_{n} \implies (iii)_n $. Let $ w\notin U_d $ but $ l(sw)=l(w)+1 $ for any $ s\in S $ with $ s\leq d $. By the Assumption \ref{ass}(i,d), we have $ dw\notin\caa $. Then by $ (ii)_n $ we have  $C_dF_w=C_{dw}=0\mod{\mc{H}_{<0}}$. Hence $ C_dT_w=-\sum_{y<_Uw}p_{y,w}C_dT_y \mod{\mc{H}_{<0}}$.  Then $ (iii)_n $ follows.   This completes a loop of  the induction. \qed

\begin{cor}\label{cor:tech}
Let $ d\in D $, $ u, v^{-1}\in U_d $.Then
\[
C_{du}=C_dF_u \mod{\mc{H}_{<\caa}}
\]
\[
C_{vd}=E_vC_d\mod{\mc{H}_{<\caa}}
\]where $ E_v=(F_{v^{-1}})^\flat $, see \ref{sec:map} for the definition of the map $ {~}^\flat $.
\end{cor}
\textit{Proof.} The first equation is just a weak version of Lemma \ref{lem:loop}(ii). Note that  $ w $ in Lemma \ref{lem:loop} is not necessarily in $ U_d $. The stronger version there are required in the proof. 

For second  equation, one only needs to apply the  anti-involution $ {~}^\flat $ to the first.\qed

\begin{cor}\label{cor:inv}Let $ d\in D $, $ u, v^{-1}\in U_d $.
\begin{itemize}
\item [(i)] There exist unique $ q_{u',u}\in\mc{A}_{\leq 0} $ with $ u'\leq _U u $ such  that\[
C_dT_u=\sum_{u'\leq _Uu}q_{u',u}C_dF_{u'}\mod{\mc{H}_{<\caa}}
\]
and $ q_{u,u}=1 $, $ q_{u',u}\in\mc{A}_{<0} $ for $ u'\neq u $.
\item [(ii)] \[
T_vC_d=\sum_{v'^{-1} \leq_U v^{-1}}q_{v'^{-1},v^{-1}}E_{v'}C_d\mod{\mc{H}_{<\caa}}.
\]
\end{itemize}
\end{cor}
\textit{Proof.} By Lemma \ref{lem:loop}(ii), we have\[
C_dT_u=C_dF_u-\sum_{u'<_Uu}p_{u',u}C_dT_{u'} \mod{\mc{H}_{<\caa}}
\]Then (i) follows immediately by induction. We get (ii)  by applying the anti-involution $ {~}^\flat $ to (i).\qed

\begin{thm}[Decomposition formula]\label{th:dec}
Keep Assumption \ref{ass}. Let $ d\in D $, $ b\in B_d $, $ u\in U_d $. We have the decomposition
\[
 C_{bdu}= E_bC_dF_u\mod{\mc{H}_{<\caa}},
 \]
where $ E_b $, $ F_u $ is given by Corollary \ref{cor:tech}. Equivalently, we have  \[C_{bd}C_{du}=h_{d,d,d}C_{bdu} \mod{\mc{H}_{<\caa}}.\]
\end{thm}
\textit{Proof.} The proof uses induction on $ l(b)+l(u) $.  By Assumption \ref{ass}(iv) and Corollary \ref{cor:inv}, we have \begin{align*}
 C_{bdu}&\equiv T_{bdu}\\
 &\edq T_bC_dT_u\\
 &\edq \sum_{b'^{-1}\leq_U b^{-1},u'\leq_U u} q_{b',b}q_{u',u}E_{b'}C_dF_{u'}.
 \end{align*}
(Recall the notation from \ref{no:equal}(iii).) By induction hypothesis, when $b'\neq b$ or $u'\neq u$, we have  $E_{b'}C_dF_{u'}\ed C_{b'du'}$  and  $(q_{b',b}q_{u',u})\in\mc{A}_{<0}$ (see Corollary \ref{cor:inv}), and hence $q_{b',b}q_{u',u}E_{b'}C_dF_{u'}\edq 0$.   So \[
 C_{bdu}\edq\sum_{b'^{-1}\leq_U b^{-1},u'\leq_U u} q_{b',b}q_{u',u}E_{b'}C_dF_{u'}\edq E_bC_dF_u.
 \]

Now by Lemma \ref{cor:tech} and Lemma \ref{lem:indepen}(ii) we have  \[C_{bd}C_{du}\ed(E_bC_d)(C_dF_u)\ed h_{d,d,d}E_bC_dF_u.\]
Then one can see that $ E_bC_dF_u $ is invariant in $\mc{H}/\mc{H}_{<c}$ under bar involution.   Therefore by Lemma \ref{lem:equal} we have\[
 C_{bdu}\ed E_bC_dF_u.
 \]
This completes the proof.\qed

\begin{cor}\label{cor:dec}
Keep Assumption \ref{ass}.
\begin{itemize}
\item [(i)] For any $ d\in D $, $ b\in U_d $, $ \Phi_{b,d}=\{bdu\mid u\in U_d \} $ is a  right cell contained in $ \caa $.
\item [(ii)]The left cells contained  in $ \caa $ are all of the form $\Theta_{b,d}=\Phi_{b,d}^{-1}  $.

\item [(iii)] If $ xy,x\in \caa $ and $ l(xy)=l(x)+l(y) $, then $ xy,x $ are in the same right cells. Similarly,  if $ xy,y\in \caa $ and $ l(xy)=l(x)+l(y) $, then $ xy,y $ are in the same left cells.

\item [(iv)] For each $ w\in\caa $, there exist unique $ d_1,d_2\in D $ and $ b_1\in B_{d_1} $, $ b_2\in B_{d_2} $ such that $ w\in \Phi_{b_1,d_1}\cap\Theta_{b_2,d_2} $.

\item [(v)] Keep the notation in (iv). There exists unique $ p_w\in \Phi_{e,d_1}\cap\Theta_{e,d_2} $ such that $ w=b_1p_wb_2^{-1} $ and $ l(w)=l(b_1)+l(p_w)+l(b_2) $. In other words, we have a union without intersection\[
\caa=\bigsqcup_{d_1,d_2\in D}B_{d_1}(\Phi_{e,d_1}\cap\Theta_{e,d_2})B_{d_2}^{-1}.
\]
\item [(vi)] Keep the notation in (v). We have decomposition
\begin{equation*}\label{eq:decbb}
C_w= E_{b_1}C_{p_w}F_{b_2^{-1}}\mod{\mc{H}_{<\caa} }.
\end{equation*}
\end{itemize}
\end{cor}
\textit{Proof.} By  Assumption \ref{ass}(iii), $ \Phi_{e,d} $ is a right cell. By decomposition formula, for any $ b\in B_d, u\in U_u $, we have $ C_{bdu}=E_bC_{du} $. Using \ref{ass}(i,d), we see immediately that $ \Phi_{b,d} $ is a right cell for any $ b\in B_d $. Then (i) and (ii) are proved.

(iii). Let $ x=bdu $. Then $ xy=bd(uy) $ with $ l(bduy)=l(b)+l(d)+l(uy) $, and hence $ uy\in U_d $. By (i), we see that $ x,xy $ are in the same right cell.

(iv). It is obvious that $ \Phi_{b_1,d_1} $ (resp. $  \Theta_{b_2,d_2} $) is the right (resp. left) cell  containing $ w $.

(v). Assume that  $ w\in \Phi_{b_1,d_1}\cap\Theta_{b_2,d_2} $. Let $ w=b_1d_1u_1 $. By (iii), the elements $ d_1u_1 $ and $ w $ are in the same left cell, i.e. $ d_1u_1 \in \Theta_{b_2,d_2} $. Write $ d_1u_1=u_2^{-1}d_2b_2^{-1} $ for some $ u_2\in U_{d_2} $. By (iii) again, the elements $ u_2^{-1}d_2 $ and $d_1u_1  $ are in the same right cell $ \Phi_{e,d_1} $. Then the element  $ p_w=:u_2^{-1}d_2 $  satisfies the conclusion in (v). 

(vi). Write $ w=b_1p_wb_2^{-1} $ as (v). By the decomposition formula, we see that $$  C_{w}\ed E_{b_1}C_{p_wb_2^{-1}}\ed E_{b_1}C_{p_w}F_{b_2^{-1}} .$$\qed

\section{Prove Assumption \ref{ass}  for affine Weyl groups of rank 2}\label{sec:assumption-holds}
\saa In this section we verify Assumption \ref{ass} for the affine Weyl groups of type $ \tilde{B}_2 $ and $ \tilde{G}_2 $. Then we can see the conclusions discussed in the last section all holds for them.
By Lemma \ref{form}, we only need to verify Assumption \ref{ass} (ii,b-c) and (iv).

\saa If $ d\in D $ is the longest element of some finite parabolic subgroup of $ W $, then Assumption \ref{ass} (ii,b-c) will automatically hold. So we only need verify the following cases for Assumption \ref{ass} (ii,b-c).

\saa The case of type $ \tilde{C}_2 $. By \ref{pc2}, there are four cases such that $ d\in D $ is not the longest element of some finite parabolic subgroup.

\begin{itemize}
 \item [(i)] When  $ c<b $ ($ a\geq c $), $101$ is an element in $ D $.
 
Since \begin{align*}
 C_1C_0C_1&=C_1T_0C_1+q^{-c}C_1C_1\\
 &=C_1T_0C_1+q^{-c}(q^{b}+q^{-b})C_1\\
 &=C_1T_0C_1+(q^{-b-c}-q^{c-b})C_1+(q^{-c+b}+q^{c-b})C_1\\
 &=C_1T_0C_1-q^{-b}\xi_cC_1+\eta_{b-c}C_1,
  \end{align*}
we can see that  $C_1T_0C_1-q^{-b}\xi_cC_1=C_1C_0C_1-\eta_{b-c}C_1$ is  bar invariant and  $\equiv T_{101} \mod{\mc{H}_{<0} }$ (see \ref{no:equal}(iv) for notations). Hence 
  \begin{align}
  C_{101}&=C_1T_0C_1-q^{-b}\xi_cC_1\label{eq:c101}\\
  &=\T{101}+q^{-b}(\T{10}+\T{01})+q^{-2b}\T{0}-q^{-b}\x{c}\T{1}-q^{-2b}\x{c}.\nonumber
  \end{align}
From \eqref{eq:c101}, we see immediately  that  $C_0C_{101}=\T{0101}\mod{\mc{H}_{<0}}$. So $ C_0C_{101}=C_{1010}$. Hence $ T_0C_{101}\ed -q^{-c}C_{101} $ (note that $ C_{1010}\in\mc{H}_{<\caa} $).  On the other hand, $T_1C_{101}=q^bC_{101}$  ($C_1C_{101}=\eta_bC_{101} $) is well-known. Then Assumption \ref{ass} (ii,b) holds in this case. From these calculations, we see also that
 \[ C_{101}C_{101}\ed (-q^{-c}\eta_b^2-q^{-b}\x{c}\eta_b)C_{101}\ed -\eta_b\eta_{b-c}C_{101}. \]

 \item [(ii)] When  $ a=c $, $ a<b $, $ 121$ is an element in $ D $.  
The calculation is very similar to  (i). One only need replace 0 by 2,  and   $c$ by $a$.
\item [(iii)] $ 212 $ is in $ D $ if and only if    $b<a$. The calculations is also similar to  (i). The results of calculations are as follows:
 \begin{align*}
 C_{212}&=C_2T_1C_2-q^{-a}\x{b}C_2,\\T_2C_{212}&=q^{a}C_{212}\\
 T_{1}C_{212}&\ed -q^{-b}C_{212}\\
 C_{212}C_{212}&\ed -\eta_a\eta_{a-b}C_{212}.
 \end{align*}
 \item [(iv)]When  $ b<c $,  $010$ is in $ D $. The calculation is similar to  (iii). One only need to change  2 to 0, and change $a$ to $c$. 
 \end{itemize}
 
 \saa {The case of type  $\tilde{G}_2$.}
By  \ref{g2} there are two cases such that  $ d\in D $ is not the longest element of a parabolic longest element.  

\begin{itemize} 
\item [(i)] When  $a>b$, 21212 is   in $ D $.
 \begin{align*}
 C_{2}C_1C_2C_1C_2&=C_2T_1T_2T_1C_2+(2q^{a-b}+q^{-a}\eta_b)C_2T_1C_2\\&+(q^{a-2b}\eta_a+\eta_a\eta_bq^{-a-b})C_2\\
 2\eta_{a-b}C_{212}& = 2\eta_{a-b} C_2T_1C_2-2\eta_{a-b}q^{-a}\x{b}C_2
 \end{align*}
So we have  \begin{align*}
 &C_{2}C_1C_2C_1C_2-2\eta_{a-b}C_{212}-(\eta_{2a-2b}+3)C_2\\
 =& C_2T_1T_2T_1C_2-q^{-a}\x{b} C_2T_1C_2+q^{-2a}(\eta_{2b}-1)C_2.
 \end{align*}
Then $ C_{21212}=C_2T_1T_2T_1C_2-q^{-a}\x{b} C_2T_1C_2+q^{-2a}(\eta_{2b}-1)C_2 $.
This implies that  $C_1C_{21212}=T_{121212}\mod{\mc{H}_{<0}}$. So $C_1C_{21212}=C_{121212}$.
And hence $T_1C_{21212}=-q^{-b}C_{21212}\mod{\mc{H}_{<c}}$. On the other hand  $T_{2}C_{21212}=q^{b}C_{21212}$ is well-known. Using these calculations we have \begin{align*}
 C_{21212}C_{21212}& \ed (q^{a-2b}\eta_a^2+q^{-a-b}\x{b}\eta_a^2+q^{-2a}(\eta_{2b}-1)\eta_a)C_{21212}\\
 &\ed \eta_a(\eta_{2a-2b}+1)C_{21212}.
 \end{align*}

 \item [(ii)] When $a<b$,  12121 is in $ D $.
The calculation is similar to (i). The results of calculations are as follows
 \begin{align*}
  C_{12121}&=C_1T_2T_1T_2C_1-q^{-b}\x{a} C_1T_2C_1+q^{-2b}(\eta_{2a}-1)C_1.\\
  T_2C_{12121}&\ed -q^{-a}C_{12121}\\
  C_{12121}C_{12121}& \ed \eta_b(\eta_{2b-2a}+1)C_{12121}.
 \end{align*}
\end{itemize}

By the above calculations, we can conclude   Assumptions \ref{ass}(ii,b-c) in present case.

\saa In remaining of this section we prove the last  assumption: \ref{ass}(iv).

The   case for lowest two-sided cell has been  known, see \cite[Cor 3.3]{xie2015lowest}. So we only   deal with non-lowest two-sided cells  with more than one element.

The computations for cases  ($ \tilde{C}_2 $,1,vii) , ($ \tilde{G}_2 $,1,i-iv), ($ \tilde{G}_2 $,2,i-ii) are more complicated than other cases. We refer cases  ($ \tilde{C}_2 $,1,vii) , ($ \tilde{G}_2 $,1,i-iv), ($ \tilde{G}_2 $,2,i-ii) as complicated cases, and refer the remaining cases as easy cases.

\saa \textbf{Easy cases}. 

In the easy cases, we actually have  conclusion  stronger than Assumption \ref{ass}(iv).
\begin{prop}\label{prop:mod}
Except the cases ($ \tilde{C}_2 $,1,vii) , ($ \tilde{G}_2 $,1,i-iv), ($ \tilde{G}_2 $,2,i-ii), we have \begin{equation}\label{nomo}
T_{b}C_dT_u=T_{bdu} \mod{\mc{H}_{<0}}
\end{equation}
Compare with Assumption \ref{ass}(iv).
\end{prop}

\begin{lemma}\label{lem:comp1}
Let $ d'<d $ with $ d\in D $, $ b\in B_d $, $ u\in U_d $ and $ \alpha\in\mc{A} $. If  $ \alpha T_{b}T_{d'}T_u\in \mc{H}_{<0} $, then for all $ b'^{-1}\leq _D b^{-1},\forall u'\leq_{D} u  $ we have $ \alpha T_{b'}T_{d'}T_{u'}\in\mc{H}_{<0} $.
\end{lemma}
\textit{Proof.} The proof is immediately from the fact that $ m_{x,y,z} $ is a polynomials in $ \xi_{L(s)},s\in S $ with positive integral coefficients, where $ T_xT_y=\sum_{z}m_{x,y,z}T_z $.\qed

The following lemma is obvious.
\begin{lemma}\label{lem:comp2}
Assume $ T_xT_y=\sum_zm_{x,y,z}T_z $. Let $ w $ be an element in $ W $ such that $ T_yT_w=T_{yw} $ and $ T_{z}T_w =T_{zw}$ for all $ z $ with $ m_{x,y,z}\neq 0 $. Then we have  $ m_{x,yw,zw}= m_{x,y,z}$. 
\end{lemma}

\begin{cor}\label{cor:inf}
let  $ d'<d \in D$, $b\in B_d, u\in U_d  $. Assume that, for any $  w \in  W $ such that $ T_{u}T_w=T_{uw} $ and $ uw\in U_d $, we have  $ T_{z}T_w =T_{zw}$ for all $ z $ with $ m_{b,d'u,z}\neq 0 $. And assume that  there exists $ \alpha\in \mc{A} $ such that $ \alpha T _bT_{d'}T_u\in\mc{H}_{<0} $. Then for all $ b'^{-1}\leq_U b^{-1} $ and $ u'\in U_d $ we have $ \alpha T_{b'}T_{d'}T_{u'}\in\mc{H}_{<0} $.
\end{cor}

\textit{Proof.} By Lemma \ref{lem:comp2}, we have $ m_{b,d'uw,zw}= m_{b,d'u,z}$. So, for all $ b'^{-1}\leq_U b^{-1} $ and $ u'\in U_d $ with $ u\leq_U u' $, we have $ \alpha T_{b'}T_{d'}T_{u'}\in\mc{H}_{<0} $. Then the corollary follows  from Lemma \ref{lem:comp1}. (Here we need a fact: for any $u, v\in U_d $ there exists  $ u' \in U_d$ such that $ u\geq_U u' $ and $ v \geq_U u' $.  For this fact, one can check easily case by case.) \qed

Now we can do some calculations case by case and then prove Assumption \ref{ass}(iv).

\noindent\textbf{Case ($ \tilde{C}_2 $,1,i)}.  $ a-c>2b $. $ B=\{e,1,01,101\} $, $ d=02 $, $ U=U(1012) $.
 \begin{align*}
 T_{101}T_0\T{1012}&=(\x{b}^2\x{c}+\x{c})\T{01012}+\x{b}\x{c}\T{0102}\\&+\x{b}^2\T{1012}+\x{b}\T{012}+\x{b}\T{102}+T_{02}\\
 \T{101}\T{2}\T{1012}&=\T{10121012}\\
 \T{101}\T{1012}&=\x{b}^2\T{210102}+\x{b}\T{0102}+\x{c}\T{1012}+\x{b}\T{12}+T_2
 \end{align*}
One can check that the conditions of Corollary \ref{cor:inf}  are satisfied. Hence 
 $ \T{b}C_dT_u =\T{bdu} \mod{\mc{H}_{<0}} ,\forall b\in B, u\in U $, i.e Proposition  \ref{prop:mod} in this case holds.

For the rest of ``easy cases",   Proposition  \ref{prop:mod} follows by similar easy calculations. The main details  are arranged in Appendix \ref{app:easy}.

\saa \textbf{The complicated cases}: ($ \tilde{C}_2 $,1,vii) , ($ \tilde{G}_2 $,1,i-iv), ($ \tilde{G}_2 $,2,i-ii). For these cases, Proposition \ref{prop:mod} does not hold; we have to return to Assumption \ref{ass}(iv). This is one of reasons why the computations are complicated.

\noindent\textbf{Cases ($ \tilde{C}_2 $,1,vii)}. $ a>c, a+c<2b $. $ B=\{e,2,12,012\} $, $ d=101 $, $ U=U(2101) $.
   \begin{align*}
   \T{012}\T{10}\T{2101}&=\T{0212010}+\x{b}\T{02121010}\\
   \T{012}\T{01}\T{2101}&=\T{010121201}\\
      \T{012}\T{0}\T{2101}&=\x{a}\T{0120101}+\x{b}\x{c}\T{1010}+\x{b}
      \T{101}+\x{c}\T{010}+\T{10}\\
      \T{012}\T{1}\T{2101}&={\x{b}\T{0121201}}+\T{021201}\\
      \T{012}\T{2101}&=\x{a}\T{012101}+\x{b}\T{1010}+\x{c}\T{01}+\T{1}
      \end{align*}    
Using Lemma \ref{lem:comp2} and the equation \eqref{eq:c101}, we have 
 \begin{align}\label{233}
\T{012}C_{101}\T{2101v}&\equiv \T{012102101v}+\T{02121010v}-q^c\T{0121201v} \qquad \forall v\in U.
    \end{align}
Note that $ 1212<02121010v $ and $ 1212< 0121201v$ and hence $ 1212<\caa $ and $ 1212< \caa$. One can check that $ \T{02121010v}\in \mc{H}_{<0}+\mc{H}_{<\caa} $ (which is easy) and $ \T{0121201v}\in \mc{H}_{<0}+\mc{H}_{<\caa} $ (using the condition $ c<a $ and $ c<b $). And hence \[
\T{012}C_{101}\T{2101v}\edq \T{012102101v} \qquad \forall v\in U.
\]Now one cannot conclude Assumption \ref{ass}(iv) like the ``easy cases", since we cannot use Lemma \ref{lem:comp1} in this case. But we only need to do more similar calculations:
\begin{align*}
\T{12}C_{101}\T{2101v}&\equiv \T{12102101v}+\T{2121010v}-q^c\T{121201v}.\\
\T{2}C_{101}\T{2101v}&\equiv \T{2102101v}.\\
\T{012}C_{101}\T{210}&\equiv \T{01210210}-q^c\T{12120}.\\
\T{012}C_{101}\T{21}&\equiv \T{0121021}-q^c\T{01212}.\\
\T{12}C_{101}\T{21}&\equiv \T{121021}-q^c\T{1212}.\label{344}
\end{align*}

Using these calculations, one can conclude Assumption \ref{ass}(iv):\[
\T{b}C_{101}\T{u}\edq \T{b101u} \qquad \forall b\in B, u\in U.
\]

The computations  for the rest ``complicated cases" are similar  and are arranged in Appendix \ref{app:complicated}.

\section{Decomposition formula and conjectures P1-P15 }

In this section, we will see that the decomposition formula will be useful for the conjectures P1-P15. One can  find our formulation  $ (P1)_{\leq\caa}$-$(P15)_{\leq\caa} $ in \ref{P1-P15}

\begin{thm}\label{th:ind}
Assume that \begin{itemize}
\item[(a)] \ref{ass} holds for the two-sided cells $ \caa $;
\item[(b)] $ \mb{a}(d)=\deg h_{d,d,d}$ for any $ d\in D$;
\item[(c)] $ (P1)_{<\caa} $, $ (P4)_{\leq\caa} $, $ (P8')_{\leq\caa} $ holds.
\end{itemize}
Then\begin{itemize}
\item[(i)] $ (P1)_{\leq\caa} $ holds.
\item[(ii)] $ \mc{D}_\caa:=\caa\cap\mc{D}$ is the set $\{bdb^{-1}\,|\,d\in D,b\in B_d \} $.
\item [(iii)] Every left cell and right cell   in $ \caa $ contains a unique element in $ \mc{D} $.
\item[(iv)]Let $ q\in\mc{D}_\caa $. Then $\gamma_{x,y,q}\neq0$ if and only if $x=y^{-1}$, $y\sim_{ {L}}q$. And in  this case we have $\gamma_{x,y,q}=n_q=\pm1$.
\end{itemize}
\end{thm}
\textit{Proof.} By Theorem \ref{th:dec}, 	
 $ C_{bd}C_{du}=h_{d,d,d}C_{bdu}\mod{\mc{H}_{<c}} $. Hence we can write \begin{equation}\label{eq:1}
 C_{bd}C_{du}=h_{d,d,d}C_{bdu}+\sum_{z<\caa}h_{bd,du,z}C_z.
 \end{equation}
 Applying the map $ \tau $ (see \ref{sec:map}) to \eqref{eq:1}, we get
 \begin{equation}\label{eq:q1}
  \delta_{bd,(du)^{-1}}=h_{d,d,d}P_{e,bdu}+\sum_{z<\caa}h_{bd,du,z}P_{e,z}\mod{\mc{A}_{<0}}
  \end{equation}
  where $ \delta_{x,y}=1 $ if  $ x=y $, otherwise $ \delta_{x,y}=0 $.
By $ (P1)_{<\caa} $, we have   $\mb{a}(z)\leq\Delta(z)$ for all   $ z<\caa $. Hence $$ \sum_{z<\caa}h_{bd,du,z}P_{e,z}=\sum_{z<\caa}\gamma_{bd,du,z^{-1}}n_z \mod{\mc{A}_{<0}}. $$
By $ (P8')_{\leq\caa} $, we have $ \gamma_{bd,du,z^{-1}} =0 $ since $ du $ and $ z $ are  in different left cells.

Therefore
\begin{equation}\label{eq:jl}
 \sum_{z<\caa}h_{bd,du,z}P_{e,z}=0 \mod{\mc{A}_{<0}}.
\end{equation}
Then equation \eqref{eq:q1} becomes
 \begin{equation}\label{eq:sdf}
 h_{d,d,d}P_{e,bdu}= \delta_{b,u^{-1}}\mod{\mc{A}_{<0}}.
 \end{equation}

So $ \Delta(bdu)=- \deg P_{e,bdu}\geq \deg h_{d,d,d}$. By the hypothesis (b) and $ (P4)_{\leq\caa} $, $\deg h_{d,d,d}=\mb{a}(d)=\mb{a}(bdu)$. Then we get (i): \[
\mb{a}(bdu)\leq\Delta(bdu).
\]
From \eqref{eq:sdf}, we can see also that  $bdu\in\mc{D}$ if and only if $b=u^{-1}$. Then  (ii) follows. 

(iii) follows from Corollary \ref{cor:dec} and (ii).

(iv). For $ y\in\caa $, we write $ \Theta_y $ the left cell containing $ y $. Denote by $ q_y $ the unique element of $ \mc{D} $ in $ \Theta_y $ (see (iii)). For any $ x,y\in\caa $, consider the equation
\begin{equation}\label{eq:jv}
  \tau(C_xC_y)=\sum_{z\in\Theta_y}h_{x,y,z}P_{e,z}+\sum_{z<\caa}h_{x,y,z}P_{e,z}.
  \end{equation}

For any $z\in \Theta_y$ with $z\neq q_y$, we have $\Delta(z)>\mb{a}(z)$ by (iii), and hence $  h_{x,y,z}P_{e,z}\in\mc{A}_{<0} $. On the other hand,  for all $ z<\caa $,  we have  $ \deg h_{x,y,z}P_{e,z}<0  $ by the same reason as \eqref{eq:jl}. Therefore, by taking constant  terms in \eqref{eq:jv}, we get
\begin{equation}\label{eq:del}
 \delta_{x,y^{-1}}=\gamma_{x,y,q_y}n_{q_y}.
\end{equation}
 Note that  $ n_{q_y}$ is a nonzero integer. 
 
 If $\gamma_{x,y,q}\neq0$, then $y\sim_{ {L}}q^{-1}=q$ by $ (P8')_{\leq\caa} $, i.e. $ q=q_y $. Then we get   $ x=y^{-1} $ by \eqref{eq:del}. Conversely, if  $ x=y^{-1} $, $y\sim_{ {L}}q^{-1}=q$, then $ q=q_y $ and  by  \eqref{eq:del} again  we get $ \gamma_{x,y,q}=\pm1=n_q $. This completes the proof of (iv).\qed
 
 \begin{thm}\label{th:commute}
Let $\mathcal{E}$ be a free $\mc{A}\otimes \mc{A}$-module with basis $\{\mathcal{E}_w\mid w\in \caa\}$. Define a left  $\mathcal{H}$-module structure on  $\mathcal{E}$ by
 \[
 C_x\mathcal{E}_w=\sum_{z\in \caa}( h_{x,w,z}\otimes 1)\mathcal{E}_z\quad\text{ for all } x\in W,w\in \caa
 \]
 and a right  $\mathcal{H}$-module structure
 \[
 \mathcal{E}_w C_y=\sum_{z\in \caa} (1\otimes h_{w,y,z})\mathcal{E}_z\quad\text{ for all }x\in W,w\in \caa.
 \]
 Assume that  Assumption \ref{ass} holds for the two-sided cells $ \caa $ and $ \mb{a}(d)=\deg h_{d,d,d}$ for any $ d\in D$. Then
 $\mathcal{E}$ is  an  $\mathcal{H}$-bimodule, i.e the left and right module structures commute. 
 \end{thm}
 
 \textit{Proof.} In this proof we will abbreviate $h'_{x,y,z}=h_{x,y,z}\otimes1\in\mc{A}\otimes\mc{A}$ and  $h''_{x,y,z}=1\otimes h_{x,y,z}\in\mc{A}\otimes\mc{A}$. 
 We first prove the following equation.  
\begin{equation}\label{eq:claim}
(C_{u^{-1}d}\mathcal{E}_{d})C_{du'}=C_{u^{-1}d}(\mathcal{E}_{d}C_{du'})\qquad \text{ for }d\in D, u, u'\in U_d.
\end{equation}

By decomposition formula, we can see $$ C_{u^{-1}d}C_{d}\ed E_{u^{-1}}C_{d}C_{d}\ed h_{d,d,d}E_{u^{-1}}C_{d}\ed h_{d,d,d}C_{u^{-1}d} $$ and hence $C_{u^{-1}d}\mathcal{E}_{d}=h'_{d,d,d}\mathcal{E}_{u^{-1}d}$. Similarly, we have  $\mc{E}_{d}C_{du'}=h''_{d,d,d}\mc{E}_{du'}$.

Since $C_{u^{-1}d}C_{du'}\ed h_{d,d,d}E_{u^{-1}}C_{d}F_{u'}$,  $E_{u^{-1}}C_{d}F_{u'} $ is bar-invariant in  $\Hm/\mc{H}_{<\caa}$. So there exists bar invariant $b_{z}\in\mc{A}$ such that $E_uC_{d}F_{u'}\ed\sum_{z\in \caa}b_zC_z$. Hence $h_{u^{-1}d,du',z}=b_zh_{d,d,d}$,  $\forall z\in \caa$. By P4 and (b) in the hypothesis, $\mathbf{a}(z)=\mathbf{a}(d)=\deg(h_{d,d,d})$. Then we have  $b_z\gamma_{d,d,d}=\gamma_{ud,du'^{-1},z^{-1}}\in \mathbb{Z}$. So $$\mathcal{E}_{u^{-1}d}C_{du'}=\frac{h''_{d,d,d}}{\gamma_{d,d,d}} \sum_{z\in \caa}\gamma_{ud,du'^{-1},z^{-1}}\mathcal{E}_z.$$ So 
 \begin{align*}
 (C_{u ^{-1}d}\mathcal{E}_{d})C_{du'}&= {h'_{d,d,d}} \mathcal{E}_{u^{-1}d}C_{du'}\\ 
 &=\frac{h'_{d,d,d}h''_{d,d,d} }{\gamma_{d,d,d}}\sum_{z\in \caa}\gamma_{u^{-1}d,du',z^{-1}}\mathcal{E}_z
 \end{align*}
 
Similarly, we have$$C_{u ^{-1}d}(\mathcal{E}_{d}C_{du'})=\frac{h'_{d,d,d}h''_{d,d,d} }{\gamma_{d,d,d}} \sum_{z\in \caa}\gamma_{u^{-1}d,du',z}\mathcal{E}_z.$$
Then equation \eqref{eq:claim} follows.
 
Let $x,y\in W$, $b\in B_d$, $u\in U_d$. By decomposition formula we have $C_{bd}\mathcal{E}_{du}=h'_{d,d,d}\mathcal{E}_{bdu}$. 
 \begin{align*}
 (C_x\mathcal{E}_{bdu})C_y&=\frac1{h'_{d,d,d}}(C_x(C_{bd}\mathcal{E}_{du}))C_y\\
 &=\frac1{h'_{d,d,d}}((C_xC_{bd})\mathcal{E}_{du})C_y\\ \intertext{(Since $(dU_d)^{-1}=\{v^{-1}d\,|\,v\in U_d \}$ is a left cell, we have $C_xC_{bd}\ed\sum_{v\in U_d}h_{x,bd,v^{-1}d}C_{v^{-1}d}$.)}
 &=\frac1{h'_{d,d,d}}(\sum_{v\in U_d}h'_{x,bd,v^{-1}d}C_{v^{-1}d}\mathcal{E}_{du})C_y \\
 &=\frac1{h'_{d,d,d}h''_{d,d,d}} (\sum_{v\in U_d}h'_{x,bd,v^{-1}d}C_{v^{-1}d}(\mathcal{E}_{d}C_{du}))C_y\\
 \intertext{by \eqref{eq:claim}}
 &=\frac1{h'_{d,d,d}h''_{d,d,d}} (\sum_{v\in U_d}h'_{x,bd,v^{-1}d}(C_{v^{-1}d}\mathcal{E}_{d})C_{du})C_y\\
 &=\frac{h'_{d,d,d}}{h'_{d,d,d}h''_{d,d,d}} (\sum_{v\in U_d}h'_{x,bd,v^{-1}d}(\mathcal{E}_{v^{-1}d}C_{du}))C_y\\
 &=\frac{h'_{d,d,d}}{h'_{d,d,d}h''_{d,d,d}} \sum_{v\in U_d}h'_{x,bd,v^{-1}d}\mathcal{E}_{v^{-1}d}(\sum_{v'\in U_d}h''_{du,y,dv'}C_{dv'})\\
 &=\frac1{h'_{d,d,d}h''_{d,d,d}} \sum_{\substack{v\in U_d\\v'\in U_d}} h'_{x,bd,v^{-1}d}h''_{du,y,dv'}(C_{v^{-1}d}\mathcal{E}_{d})C_{dv'}
 \end{align*}
 
Similarly, \[C_x(\mathcal{E}_{bdu}C_y)=\frac1{h'_{d,d,d}h''_{d,d,d}} \sum_{\substack{v\in U_d\\v'\in U_d}} h'_{x,bd,v^{-1}d}h''_{du,y,dv'}C_{v^{-1}d}(\mathcal{E}_{d}C_{dv'}).\] By \eqref{eq:claim} again, we have $(C_x\mathcal{E}_{bdu})C_{y}=C_x(\mathcal{E}_{bdu}C_y)$. Then the theorem follows.\qed

 \begin{thm}\label{th:conjrk2}
Let  $ W $ be an affine Weyl group of type $ \tilde{C}_2 $ or  $ \tilde{G}_2 $. If  P4, P8' hold and  $ \mb{a}(d)=\deg h_{d,d,d} $ for all $ d\in D $, then P1-P15 hold.
 \end{thm}
 \textit{Proof.} First P14 is obvious from the partition of two-sided cell of Guilhot \cite{guilhot2010rank2,guilhot0810.5165}, see Section \ref{p-c}.
  Recall that Assumption \ref{ass} holds for the  affine Weyl group of type $ \tilde{C}_2 $ or  $ \tilde{G}_2 $, see Section \ref{sec:assumption-holds}. Then we can apply Theorem \ref{th:commute} and we get immediately P15.
 
 By induction,  using Theorem  \ref{th:ind}  we obtain  P1, P2, P3, P5, P6, P13 immediately.
 
 By P4 and  the assumption that $ \mb{a}(d)=\deg h_{d,d,d} $ for all $ d\in D $, we can determine the $ \mb{a} $-values on each two-sided cell. By  explicit computations of $ \deg h_{d,d,d} $ (see Section \ref{sec:assumption-holds}), we can conclude that $ \mb{a}(\caa')<\mb{a}(\caa) $ for $ \caa'<\caa $. Then we get P11.
 
 Using  P11 and \cite[Rem.6.7]{guilhot2010rank2}, we get P9, P10. And   P3, P4, P8' imply  P12, see \cite[\S14.12]{lusztig2003hecke}; P2, P3, P4, P5 imply P7, see\cite[\S14.7]{lusztig2003hecke}; P7 implies P8, see \cite[\S14.8]{lusztig2003hecke}. This completes the proof.\qed
 
 \begin{remark}
 It is an interesting question to compute the $ \mb{a} $-function  for the affine Weyl groups of type $ \tilde{C}_2 $ or  $ \tilde{G}_2 $. If the $ \mb{a} $-values are determined, then the assumptions in Theorem \ref{th:conjrk2} will be clear except P8'. 
 \end{remark}

\appendix
\section{Computations for easy cases.}\label{app:easy}

\noindent\textbf{Case ($\tilde{C}_2$,1,ii)}. $ 0<a-c<2b $. $ B=\{e,,2,12,012 \} $, $ d=1010 $, $ U=U(2101) $.

 \begin{align*}
\T{012}\T{101}\T{2101}&=\T{0121012101}\\
\T{012}\T{010}\T{2101}&=\T{1010212010}\\
\T{012}\T{10}\T{2101}&=\x{b}\T{01212010}+\T{0212010}\\
\T{012}\T{01}\T{2101}&=\T{01021201}\\
\T{012}\T{0}\T{2101}&={\x{a}\T{0102101}}+\x{b}\x{c}\T{1010}+\x{b}\T{101}+\x{c}\T{010}+\T{10}\\
\T{012}\T{1}\T{2101}&=\x{b}\T{0121201}+\T{021021}\\         
\T{012}\T{2101}&={\x{a}\T{012101}}+\x{b}
\T{1010}+\x{c}\T{01}+T_1.
 \end{align*}
 
 \noindent\textbf{Case ($\tilde{C}_2$,1,iii)}. $ |a-c|<b, a+c>b $. $ B=\{e,1,21,01\} $, $ d=02 $, $ U=U(102) $.
  \begin{align*}
  \T{21}\T{2}\T{102}&={\x{a}\T{{1212}0}}+\T{1210}\\
  \T{21}T_0\T{102}&=\T{210102}\\
  \T{21}\T{102}&=\x{b}\T{2102}+\x{a}\T{02}+T_0.
  \end{align*}
  
\noindent\textbf{Case ($\tilde{C}_2$,1,iv)}. $ a+c<b $. $ B=\{e,0,2,02\}$, $ d=1 $, $ U=U(021) $.
  \begin{align*}
   \T{02}\T{021}&= \x{a}\x{c}\T{021}+\x{c}\T{01}+\x{a}\T{21}+\T{1}\\
  \end{align*}

\noindent\textbf{Case ($\tilde{C}_2$,1,v)}. $ a-c>b $. $ B=\{e,0,10,010\} $, $ d=212 $, $ U=U(012) $.
  
  \begin{align*}
  \T{010}\T{12}\T{012012}&=\x{c}\T{1010212012}+\T{101212012}\\
  \T{010}\T{21}\T{012012}&=\x{c}\T{0121010212}+\T{012101212}\\
  \T{010}\T{1}\T{012012}&=(\x{b}\x{c}^{2}+\x{b}\x{c})\T{101212}+\x{c}^{2}\T{010212}+
  \x{c}\T{10212}\\&+\x{b}\T{0101212}+\x{c}\T{01212}+\T{1212}\\
  \T{010}\T{2}\T{012012}&=\x{c}\T{010210212}+\T{01212012}\\
\T{010}\T{012012}&=\x{c}^2+\x{c}\T{101212}+\x{b}\T{012012}+\x{c}\T{2012}+\T{212}
\end{align*}
  
\noindent\textbf{Case ($\tilde{C}_2$,1,vi)}. $ a>c, a+c>2b $. $ B=\{e,1,01,101\} $, $ d=2 $, $ U=U(1012) $. \begin{align*}
    \T{101}\T{1012}&= \x{b}^2\T{10102}+\x{b}\T{0120}+\x{c}\T{1012}+{\x{b}\T{12}}+T_2
  \end{align*}
  
  \noindent\textbf{Case ($\tilde{C}_2$,2,i-vi)}. The verification is very easy.  
  
  \noindent\textbf{Case ($\tilde{C}_2$,3,i)}. $ a-c=2b $. $ B_1=\{e,1,01 \}, d_1=02, U_1=U(1012) $, $B_2=\{e\}$.
   
   \begin{align*}
   \T{01}T_0\T{101}&=\x{b}\x{c}\T{1010}+\x{b}\T{101}+\x{c}\T{010}+\T{10}\\
   \T{01}T_2\T{101}&=\T{012101}\\
   \T{01}\T{101}&=\x{b}\T{1010}+\x{c}\T{01}+\T{1}.
   \end{align*}
   
\noindent\textbf{Case ($\tilde{C}_2$,3,ii)}. $a+c=b$. $ B_1=\{e,0,2 \}, d_1=1, U_1=U(021) $, $B_2=\{e\}$.
  \begin{align*}
  T_0T_{02}&=\x{c}T_{02}+T_2
  \end{align*}

\noindent\textbf{Case ($\tilde{C}_2$,3,iii)}.  $a-c=b$. $ B_1=\{e,1,01 \}, d_1=02, U_1=U(102) $, $B_2=\{e\}$.  
  \begin{align*}
  \T{01}C_{02}\T{102102}&=\T{0102102102}+q^{-c}\T{012102102}+q^{-a}\x{c}\T{10102120}+q^{-a}\T{1012120}\\
  &+q^{-a-c}\x{b}\T{0102120}+q^{-a-c}\x{c}T_{02102}+q^{-a-c}\T{2102}
  \end{align*}
  
\noindent\textbf{Case ($\tilde{C}_2$,3,iv)}. $ a>c, a+c=2b $. $ B_1=\{e,1,01 \}, d_1=2, U_1=U(1012) $, $B_2=\{e\}$.
    \begin{align*}
    \T{01}\T{1012}&=\x{b}\T{01012}+\x{c}\T{012}+\T{12}.
    \end{align*}
    
\noindent\textbf{Case ($\tilde{C}_2$,3,v)}. $ a=c, a>b $. $ B_1=B_2=\{e,1\} $, $ d_1=2 $, $ U_1=U(1012) $, $ d_2=0 $, $ U_2=U(1210) $.
     \begin{align*}
     T_1C_2T_{1012}&=\T{121012}+q^{-a}T_1^2\T{012}.
     \end{align*}
     
\noindent\textbf{Case ($\tilde{C}_2$,3,vi)}. $ a=c, a<b $. $ B_1=\{e,2\} $, $ d_1=101$, $ U_1=U(2101) $, $ B_1=\{e,0\} $, $ d_2=121 $, $ U_2=U(0121) $.
  \begin{align*}
  \T{2}\T{10}\T{2101}&=\T{2120101}\\
   \T{2}\T{01}\T{2101}&=\T{0212101}\\
   \T{2}\T{0}\T{2101}&=\T{_2}^2\T{1010}\\
   \T{2}\T{1}\T{2101}&=\T{121201}\\
   \T{2}\T{2101}&=\T{2}^2\T{101}
  \end{align*}
  
  \noindent\textbf{Case ($\tilde{C}_2$,4,i-ii)}. $ \forall d\in D $, $ B_d=\{e\} $  and hence there is nothing need to verify.  
  
  \noindent\textbf{Case ($\tilde{C}_2$,4,iii)}. We only need to verify  that $ T_0C_1\T{210} \equiv  \T{01210} $, which is easy. 
  
\noindent\textbf{Case ($\tilde{G}_2$,1,v)}. $ 2a>3b,a<2b $. $ c=\{e,1,21,121,0121\}02\{e,1,12,121,1210\}$. 
  \begin{align*}
  \T{0121}\T{0}\T{1210}&=\T{012101210}\\
   \T{0121}\T{2}\T{1210}&=\x{b}\T{01212120}+\T{0212120}\\
    \T{0121}\T{1210}&=\x{b}\T{0121210}+\x{a}\T{01210}+\x{b}\T{010}+\x{b}\T{0}+1
 \end{align*}

\noindent\textbf{Case ($\tilde{G}_2$,1,vi)}. $ a<b $. $   c=\{e,0\}12121\{e,0\}$.
 \[
 \T{0 }C_{12121}\T{ 0}\equiv\T{0121210}+(q^{-3b}(\eta_{2a}-1)-q^{-3b}\x{a}\T{2})T_2^2\equiv \T{0121210}.
 \]

\noindent\textbf{Case ($\tilde{G}_2$,3,i)}. $ \forall d\in D $, $ B_d=\{e\} $ and hence there is nothing need to verify. 

\noindent\textbf{Case ($\tilde{G}_2$,3,ii)}. We only need to compute $T_2C_1\T{212}\equiv\T{21212} $.

\noindent\textbf{Case ($\tilde{G}_2$,3,iii)}. $ \forall d\in D $, $ B_d=\{e\} $ and hence there is nothing need to verify. 

\noindent\textbf{Case ($\tilde{G}_2$,3,iv)}.  $ \T{01}C_2\T{10}\equiv \T{01210} $.

\section{Computations for complicated cases}\label{app:complicated}

\noindent\textbf{Cases ($ \tilde{G}_2 $,1,i)}.  $ 2a>3b $. $ B=\{e,0,10,210,1210,01210 \} $, $ d=21212 $, $ U=U(01212) $.

\begin{align*}
T_{01210}T_{2121}T_{0121201212}&=T_{0121021210121201212}\\
T_{01210}T_{1212}T_{0121201212}&=T_{0121012120121201212}\\
T_{01210}T_{212}T_{0121201212}&=T_{012102120121201212}\\
T_{01210}\T{121}\T{0121201212}&=\x{b}\T{01201210121201212}+{\x{a}\x{b}\T{12012121012121}}\\
&+\x{a}\T{1210121201212}+\x{b}\T{1012121012121}+\T{102121012121}\\
\T{01210}\T{12}\T{0121201212}&=\x{b}\T{0121201212101212}+\underline{\x{a}\x{b}\T{1021212101212}}\\
&+\x{a}\T{102121201212}+\x{b}\T{101212101212}+\T{10212101212}\\
\T{01210}\T{21}\T{0121201212}&=\x{b}\T{0121210121201212}+\underline{\x{a}\x{b}\T{0212121012121}}\\
&+\x{a}\T{021212012121}+\x{b}\T{012121021212}+\T{02121021212}\\
\T{01210}\T{2}\T{0121201212}&=\x{b}\T{012102121201212}+\underline{\x{a}\x{b}\T{021212101212}}\\
&+\underline{\x{a}\T{02121201212}}+\x{b}\T{01212101212}+\T{0212101212}\\
\T{01210}\T{1}\T{0121201212}&=\x{b}^{2}\T{01210121201212}+\x{b}\T{0120121201212}\\
&+\x{b}\T{0121021201212}+\x{a}\T{01021201212}+\x{b}^2\T{01021212}\\
&+\x{b}\T{1021212}+\x{b}\T{0121212}+\T{212121}\\
\T{01210}\T{0121201212}&=\x{b}T_{01210121201212}+\x{b}\T{012121201212}\\
&+\x{a}\T{0121201212}+\x{b}\T{01201212}+\x{b}\T{021212}+\T{21212}.
\end{align*}
Since $ C_{21212}=C_2T_1T_2T_1C_2-q^{-a}\x{b} C_2T_1C_2+q^{-2a}(\eta_{2b}-1)C_2 $, we have 
\begin{multline*}\label{1212}
T_{01210}C_{21212}\T{0121201212u}\equiv T_{01210 21212 0121201212u}+q^{2b-a}\T{02121201212u}\\+q^{3b-a}\T{012121201212u}
-q^{2b-a}(T_{1012121201212u}+\T{0121212012121u})
\end{multline*}
Further computations shows that 
\begin{align*}
T_{01210}C_{21212}\T{012120v}&\equiv T_{01210 21212 012120v}+q^{2b-a}\T{0212120v}\nonumber\\
&+q^{3b-a}\T{01212120v}-q^{2b-a}(T_{101212120v}+\T{012121012v})\\
T_{01210}C_{21212}\T{01212 }&\equiv T_{01210 21212 01212 }+q^{2b-a}\T{021212 }\nonumber\\
&+q^{3b-a}\T{0121212 }-q^{2b-a}T_{10121212 }\\
\T{01210}C_{21212}\T{01210}&\equiv \T{012102121201210}\\
\T{1210}C_{21212}\T{0121201212u}&\equiv\T{1210212120121201212u}+q^{2b-a}\T{2121201212u}\nonumber \\&+q^{3b-a}\T{21212101212u}-q^{2b-a}\T{121212012121u}.\\
\T{1210}C_{21212}\T{01212w}&\equiv\T{12102121201212w}+q^{2b-a}\T{21212w}+q^{3b-a}\T{212121w}\\
\T{210}C_{21212}\T{01212z}&\equiv\T{2102121201212v}+q^{2b-a}\T{121212z}\\
\T{10}C_{21212}\T{0121201212u}&\equiv\T{10212120121201212u}
\end{align*}
where $ v\in\{e,1,12,121\}$, $ w\in \{e,0,01,012,0121 \} $, $ z\in \{e,0,01,012,0121,01212\} $. Using condition $ 2a>3b $, one can check that for any $u\in U$,
\begin{align*}
q^{2b-a}T_{121212u}&\equiv q^{2b-a}C_{121212}T_u\\ 	
q^{2b-a}\T{21212u}+q^{3b-a}\T{121212u}&\equiv q^{2b-a}C_{121212}\T{u}.\\
q^{2b-a}T_{0121212u}&\equiv q^{2b-a}T_0C_{121212}T_u\\ q^{2b-a}T_{10121212u}&\equiv q^{2b-a}T_{10}C_{121212}T_u\\ 		
q^{2b-a}\T{021212u}+q^{3b-a}\T{0121212u}&\equiv q^{3b-a}T_0C_{121212}\T{u}.\\
(\text{When } 2b>a  ) ~q^{2b-a} \T{012121012v}&\equiv q^{2b-a}\T{01212}C_{101} \T{2v} 
\end{align*}

 Then we can conclude that   $ \T{b}C_dT_u \edq\T{bdu}$, for any $b\in B, u\in U $.
 
 \noindent\textbf{Case ($ \tilde{G}_2 $,1,ii)}.  $ 2a<3b $.  $ B=\{e, 1, 21, 121, 0121, 2121 \} $, $ d=02 $, $ U=U(12120) $.
 
 \begin{align*}
 T_{0121}\T{0}\T{12120}&=\T{0121012120}\\
 T_{0121}\T{2}\T{12120}&=\underline{\x{a}\x{b}\T{01212120}+\x{a}\T{0212120}+\x{b}\T{0121210}}+\T{021210}\\
 T_{0121}\T{12120}&=T_2+\x{b}\T{02}+\x{b}\T{1012}+\x{a}\T{012120}+\x{b}\T{01212120}
 \end{align*}
 
Then for all $ u\in U $ we have \[
 T_{0121}C_{02}\T{12120u}\equiv\T{01210212102u}+q^{a-b}\T{0212120u}+q^{a}\T{01212120u}+\T{0121210u}.
 \]
Further computations shows that \begin{align*}
T_{0121}C_{02}\T{1212v}&\equiv\T{0121021212v}+q^{a-b}\T{021212v}+q^{a}\T{0121212v}+\T{012121v}.\\
T_{0121}C_{02}\T{1210}&\equiv \T{0121021210}.\\
T_{121}C_{02}\T{1212v}&\equiv\T{121021212v}+q^{a-b}\T{21212v}+q^{a}\T{121212v}+\T{12121v}.\\
T_{21}C_{02}\T{1212v}&\equiv\T{21021212v}+q^{a-b}\T{121212v}.\\
T_{1}C_{02}\T{1212v}&\equiv\T{1021212v}.
\end{align*}where $ v $ is the elements such that $ 1212v\in U $.
By 
\begin{align*}
q^aC_{0121212v}&\equiv q^{a-b}\T{021212v}+q^{a}\T{0121212v}+\T{012121v},\\
q^aC_{121212v}&\equiv q^{a-b}\T{21212v}+q^{a}\T{121212v}+\T{12121v},\\
 q^{a-b}C_{121212v}&\equiv q^{a-b}\T{121212v},
\end{align*} for all $ b=e,0121,121,21,1 $, $ \forall u\in U $, we have $ \T{b}C_dT_u \edq\T{bdu} $ (note that $ 121212<\caa $). It remains to verify Assumption \ref{ass}(iv) for  $ b=2121 $.
Since \begin{align*}
T_{2121}\T{0}\T{12120}&=\T{2121012120}\\
T_{2121}\T{2}\T{12120}&=\T{1210}+\underline{\x{a}\T{12120}+\x{a}T_{21210}}\\&+\underline{\x{a}\x{b}\T{121210}+\x{a}^2\T{212120}+\x{a}^2\x{b}\T{1212120}}\\
T_{2121}\T{12120}&=\T{0}+\x{a}\T{02}+\x{b}\T{2120}\\&+\x{b}\T{121210}+\x{a}\T{212120}+\underline{\x{a}\x{b}\T{1212120}}.
\end{align*} 
we have for all  $ u\in U $, \begin{align*}
T_{2121}C_{02}\T{12120u}&\equiv\T{21210212120u}+\T{1212120u}+(q^{2a}-2-q^{2a-2b})\T{1212120u}\\&+q^{a}\T{121210u}+q^{2a-b}\T{212120u}+q^{a-b}\T{21210}+q^{a-b}\T{12120u}.
\end{align*}
Further computations show that  \begin{align*}
T_{2121}C_{02}\T{1212v}&\equiv\T{2121021212v}+\T{121212v}+(q^{2a}-2-q^{2a-2b})\T{121212v}\\&+q^{a}\T{12121v}+q^{2a-b}\T{21212v}+q^{a-b}\T{2121v}+q^{a-b}\T{1212v}.
\end{align*}where $ v $ is  an element such that  $ 1212v\in U $.
Now by  \begin{align*}
q^{2a}C_{121212}T_v&\equiv q^{2a}\T{121212v} +q^{a}\T{12121v}+q^{2a-b}\T{21212v}\\&+q^{a-b}\T{2121v}+q^{a-b}\T{1212v}\\
(2+q^{2a-2b})C_{121212}T_v&\equiv (2+q^{2a-2b})T_{121212v},
\end{align*}
we can see Assumption \ref{ass}(iv) holds for   $ b=2121, u\geq_D 1212 $. 

\noindent\textbf{Case ($ \tilde{G}_2 $,1,iii)}. $ a>2b $. $ B=\{e, 1, 01, 21, 121,0121\} $, $ d=02 $, $ U=U(102) $.
\begin{align*}
\T{0121}T_0\T{102102}&=\T{120121210}+\x{b}\T{1012121012}\\
\T{0121}T_2\T{102102}&=\T{201212012}+\underline{\x{b}\T{0212121012}}\\
\T{0121}\T{102102}&=\T{12}+\x{b}\T{120}+\x{b}\T{012}+\x{b}^2\T{1012}\\&+\x{a}\T{1012120}+\x{b}\T{012121012}.
\end{align*}
Then $ T_{0121}C_{02}\T{102102u}\equiv\T{012102102102u}+\T{0121212012u},~\forall~  u\in U $.
Further computations show that 
\begin{align*}
 T_{0121}C_{02}\T{10210v}&\equiv\T{01210210210v}+\T{012121201v}\\
  T_{0121}C_{02}\T{10212}&\equiv\T{01210210212}\\
  T_{121}C_{02}\T{10210v}&\equiv\T{1210210210v}+\T{12121201v}\\
  T_{21}C_{02}\T{10210v}&\equiv\T{210210210v}\\
  T_{01}C_{02}\T{10210v}&\equiv\T{010210210v}
\end{align*}
So $ \T{b}C_dT_u \edq\T{bdu},~\forall~ b\in B, u\in U $. 

\noindent\textbf{Case($ \tilde{G}_2 $,1,iv)}. $ a<2b $. $ B=\{e,2,12,212,1212,01212\} $. $ d=101 $. $ U=U(210) $. 
 \begin{align*}
\T{01212}\T{10}\T{210210}&=\T{02121201210}+\underline{\x{b}\T{012121201210}}\\
\T{01212}\T{01}\T{210210}&=\T{01212 01 210210}\\
\T{01212}\T{1}\T{210210}&=\underline{\x{a}(\x{b}^2+1)\T{012121201}+\x{a}\x{b}\T{02121201}}+\\&\underline{\x{b}^2\T{01212101}}+\x{b}\T{0212101}+\x{b}\T{0121201}+\T{021201}\\
\T{01212}\T{0}\T{210210}&=\T{10121210}+\x{b}\T{101212101}+\x{a}\T{01212101210}\\
\T{01212}\T{210210}&=T_1+\x{b}T_{10}+\x{b}\T{01}+\x{b}^2\T{101}+\x{a}\T{101210}\\&+\x{a}\T{20121201}+\x{b}\T{01212101}+\x{a}\x{b}\T{201212101}
 \end{align*}
Then for all $ u\in U $, \begin{align*}
 \T{01212}C_{101}\T{210210u}&\equiv \T{01212101210210u}+\T{012121201210u}+\T{01212101u}\\&+q^{a-b}\T{02121201u}+q^a\T{012121201u}.
 \end{align*}
Further computations show that
\begin{align*}
\T{01212}C_{101}\T{2102z}&\equiv \T{012121012102z}+q^{a-b}\T{20121210z}+\T{0212121012z} (z\in \{e,1\})\\
\T{01212}C_{101}\T{210}&\equiv \T{01212101210}+\T{021212101}\\
 \T{1212}C_{101}\T{210210u}&\equiv \T{1212101210210u}+\T{12121201210u}+\T{1212101u}\\&+q^{a-b}\T{2121201u}+q^a\T{12121201u}\\
  \T{1212}C_{101}\T{2102w}&\equiv \T{12121012102w}+q^{a-b}\T{1212120w} (w\in \{e,1\})\\
\T{212}C_{101}\T{210210u}&\equiv \T{212101210210u}+q^{a-b}\T{12121201u} (u\in U)\\
\T{212}C_{101}\T{21021}&\equiv \T{21210121021}\\
\T{12}C_{101}\T{210210u}&\equiv \T{12101210210u}(u\in U).
\end{align*}
Since \begin{align*}
q^aC_{012121201u}&\equiv q^{a-b}\T{02121201u}+q^a\T{012121201u},\\
q^aC_{12121201u}&\equiv q^{a-b}\T{2121201u}+q^a\T{12121201u},\\
q^{a-b}C_{121212}&\equiv q^{a-b}\T{121212},
\end{align*}
we have   $ \T{b}C_dT_u \edq\T{bdu},\forall ~b\in B, d\in U $.

\noindent\textbf{Case($ \tilde{G}_2 $,2,i)}. It is part of  computations of Case ($ \tilde{G}_2 $,1,ii).The detail is omitted. 

\noindent\textbf{Case($ \tilde{G}_2 $,2,ii)}. It is part of  computations of Case ($ \tilde{G}_2 $,1,iii). The detail is omitted.

\end{document}